\documentclass[11pt,a4paper,reqno,twoside,psamsfonts]{amsart}
\usepackage[latin1]{inputenc}
\usepackage[T1]{fontenc}
\usepackage{textcomp}
\usepackage[british]{babel}
\usepackage{amscd}
\usepackage{pstricks,pst-plot}
\usepackage{xspace}
\usepackage{mathrsfs} 
\usepackage{amsmath}
\usepackage{amssymb}
\usepackage{amsthm}
\usepackage{amsfonts}
\newcommand{\CC}{{\ensuremath{\mathord{\mathbb{C}}}}\xspace}
\newcommand{\RR}{{\ensuremath{\mathord{\mathbb{R}}}}\xspace}
\newcommand{\QQ}{{\ensuremath{\mathord{\mathbb{Q}}}}\xspace}
\newcommand{\DD}{{\ensuremath{\mathord{\mathbb{D\mkern0.4\thinmuskip}}}}\xspace}

\newcommand{\NN}{{\ensuremath{\mathord{\mathbb{N}}}}\xspace}

\newcommand{\II}{{\ensuremath{\mathord{\mathbb{I}}}}\xspace}
\newcommand{\ZZ}{{\ensuremath{\mathord{\mathbb{Z}}}}\xspace}
\DeclareMathAlphabet{\Bmi}{OT1}{cmm}{b}{it}
\newcommand{\EPS}{\ensuremath{\varepsilon}\xspace}
\newcommand{\IM}{\operatorname{Im}}
\newcommand{\RE}{\operatorname{Re}}
\newcommand{\supp}{\operatorname{supp}}

\newcommand{\dbar}{\ensuremath{\mathord{\overline{\partial}}}\xspace}

\newcommand{\FP}{\operatorname{f.p.}}
\newcommand{\MTPI}{\frac{-1}{2\pi i}}
\newcommand{\ABS}[1]{{\ensuremath{\mathord{\left|#1\right|}}}}
\newcommand{\NORM}[1]{{\mathord{\left\|#1\right\|}}}
\newcommand{\DEF}{\mathbin{\smash[t]{\overset{\scriptscriptstyle\mathrm{def}}{=}}}}
\newcommand{\CTENS}[1][]{\mathbin{\widehat{\otimes}_{#1}}}
\newcommand{\PMI}{{±\infty}}
\newcommand{\lirarrow}{\stackrel{\sim}{\longrightarrow}} 
\newcommand{\SPROD}[2]{
            \ensuremath{\mathord{\left\langle#1\,\mathord{,}\,#2\right\rangle}}\xspace}

\newcommand{\dd}{\ensuremath{\mathrm{d}}} 
\newcommand{\ee}{\ensuremath{\mathrm{e}}} 
\newcommand{\ie}{, {\it i.e.},\xspace}
\newcommand{\eg}{, {\it e.g.},\xspace}
\newcommand{\IX}[1]{{\ensuremath{\mathord{\mathbf{#1}}}}} 
\newcommand{\AIX}[1]{\ABS{\mathbf{#1}}}                   
\newcommand{\GLQQ}{\char"12\kern .08em} 
\newcommand{\GRQQ}{\kern .08em\char"10\xspace}
\def\SetConstructor#1#2#3#4{%
  \def\test{#4}\ifx\test\empty{\ensuremath{\mathord{#1}_{#2}^{#3}}\xspace}%
  \else{\ensuremath{\mathord{#1}_{#2}^{#3}(#4)}\xspace}\fi}
\newcommand{\SIDX}[3]{\SetConstructor{\mathscr{S}}{#1}{#2}{#3}}
\newcommand{\SD}[2][]{\SIDX{#2}{}{#1}}                         
\theoremstyle{plain}
\newtheorem{theo}{Theorem}[section]
\newtheorem*{theo*}{Theorem}
\newtheorem{defitheo}[theo]{Definition and Theorem}
\newtheorem{prop}[theo]{Proposition}
\newtheorem{coro}[theo]{Corollary}
\newtheorem{lemm}[theo]{Lemma}
\theoremstyle{definition}
\newtheorem{defi}[theo]{Definition}
\newtheorem{defirem}[theo]{Definition and Remark}
\newtheorem{rem}[theo]{Remark}

\theoremstyle{remark}
\newtheorem*{note}{Note}

\newtheorem{nxmpl}[]{Example}
\newcommand{\labelT}[1]{\label{theo:#1}}

\newcommand{\labelL}[1]{\label{lemm:#1}}

\newcommand{\refT}[1]{Theorem~\ref{theo:#1}}

\newcommand{\refL}[1]{Lemma~\ref{lemm:#1}}

\usepackage{hyperref}
\usepackage{mathptmx,courier}
\newif\ifpdf
\ifx\pdfoutput\undefined
  \pdffalse
\else
  \pdftrue
\fi
\ifpdf 
  \hypersetup{%
    pdftitle= {Asymptotic Hyperfunctions, Tempered Hyperfunctions, and Asymptotic Expansions},
    pdfauthor=Andreas U. Schmidt,
    pdfsubject= {2000 Mathematical subject classification: Primary: 46F15; Secondary: 46F20, 30E15},
    pdfkeywords= {Hyperfunctions, asymptotic expansions, moments, Radon transformation} }
\else 
  \providecommand{\href}[2]{#2}
\fi
\newgray{mid1}{0.875}
\newgray{mid2}{0.4}
\newpsobject{showgrid}{psgrid}{subgriddiv=1,griddots=10,gridlabels=6pt}
\psset{framesep=2pt}
\psset{dotsep=2pt}
\newcommand{\QIDX}[3][]{\SetConstructor{\mathscr{Q}}{#2}{#1}{#3}}
\newcommand{\FQIDX}[2]{\SetConstructor{\widehat{\mathscr{Q}}}{#1}{}{#2}}
\newcommand{\BIDX}[2]{\SetConstructor{\mathscr{B}}{#1}{}{#2}}
\newcommand{\PIDX}[3]{\SetConstructor{\mathscr{P}}{#1}{#2}{#3}}
\newcommand{\DIDX}[3]{\SetConstructor{\mathscr{D}}{#1}{#2}{#3}}
\newcommand{\KIDX}[3]{\SetConstructor{\mathscr{K}}{#1}{#2}{#3}}

\newcommand{\HOIDX}[3]{\SetConstructor{\mathcal{O}}{#1}{#2}{#3}}
\newcommand{\HOIDXOL}[3]{%
\SetConstructor{\overline{\mathcal{O}}}{#1}{#2}{#3}}

\newcommand{\SHOPI}{\SetConstructor{\mathcal{O}}{\infty}{}{}}
\newcommand{\SHOMI}{\SetConstructor{\mathcal{O}}{-\infty}{}{}}
\newcommand{\SHQ}[1][\gamma]{\QIDX{#1}{}}

\newcommand{\SHQPI}{\QIDX{\infty}{}}
\newcommand{\SHQMI}{\QIDX{-\infty}{}}
\newcommand{\SHP}[1][\gamma]{\PIDX{#1}{}{}}
\newcommand{\SHPPI}{\PIDX{\infty}{}{}}
\newcommand{\SHPMI}{\PIDX{-\infty}{}{}}
\newcommand{\SQ}[2][\gamma]{\QIDX{#1}{#2}}

\newcommand{\SQPI}[1]{\QIDX{\infty}{#1}}
\newcommand{\SQMI}[1]{\QIDX{-\infty}{#1}}
\newcommand{\SB}[2][]{\BIDX{#1}{#2}}

\newcommand{\SPPI}[1]{\PIDX{\infty}{}{#1}}
\newcommand{\SMODPPI}[1]{\SetConstructor{\widetilde{\mathscr{P}}}{\infty}{}{#1}}
\newcommand{\SMODQMI}[1]{\SetConstructor{\widetilde{\mathscr{Q}}}{-\infty}{}{#1}}
\newcommand{\SPMI}[1]{\PIDX{-\infty}{}{#1}}
\newcommand{\SK}[2][\gamma]{\KIDX{#1}{}{#2}}
\newcommand{\SSS}[2][]{\SIDX{#1}{}{#2}}

\renewcommand{\SD}[2][]{\DIDX{}{#1}{#2}}

\newcommand{\SDDFO}[2][\beta]{\DIDX{#1}{\mathrm{F}\prime}{#2}}
\newcommand{\SSD}[2][]{\SIDX{}{#1\prime}{#2}}
\newcommand{\SKD}[2][]{\KIDX{}{#1\prime}{#2}}

\newcommand{\HOIE}[1]{\SetConstructor{\mathcal{O}}{}{\ast}{#1}}
\newcommand{\CIIED}[1]{\SetConstructor{C}{}{\infty\ast\prime}{#1}}
\newcommand{\CIE}[2][\infty]{\SetConstructor{C}{}{#1\ast}{#2}}
\newcommand{\CIIE}[1]{\CIE{#1}}

\newcommand{\HOBO}[2][\gamma]{\HOIDXOL{#1}{B}{#2}}

\newcommand{\HOPI}[1]{\HOIDX{\infty}{}{#1}}
\newcommand{\HOMI}[1]{\HOIDX{-\infty}{}{#1}}
\newcommand{\HOPMI}[1]{\HOIDX{±\infty}{}{#1}}

\newcommand{\HOPMID}[1]{\HOIDX{±\infty}{\prime}{#1}}
\newcommand{\PPMID}[1]{\PIDX{±\infty}{\prime}{#1}}
\newcommand{\PPMI}[1]{\PIDX{±\infty}{}{#1}}
\newcommand{\QPMI}[1]{\QIDX{±\infty}{#1}}
\newcommand{\NO}[3]{%
\ensuremath{\mathord{\overline{\NORM{#1}}_{{#2},{#3}}}}\xspace}

\newcommand{\NODX}[3]{%
\ensuremath{\mathord{\NORM{#1}}_{{#2},{#3}}}\xspace}

\newcommand{\RELCO}[4]{\SetConstructor{\mathrm{H}}{#2}{#1}{#3;#4}}

\newcommand{\MOM}[2][]{\SetConstructor{\mu}{}{#2}{#1}}

\newcommand{\HELMOM}[3][\omega]{\SetConstructor{p}{\RAD #2}{#3}{#1}}

\newcommand{\SUM}[2][]{\SetConstructor{S}{#1}{#2}{}}
\newcommand{\REM}[2][]{\SetConstructor{R}{#1}{#2}{}}
\newcommand{\CES}[1][]{%
  \def\test{#1}\ifx\test\empty{\ensuremath{\mathord{(\mathrm{C})}}\xspace}%
  \else{\ensuremath{\mathord{(\mathrm{C},#1)}}\xspace}\fi}
\newcommand{\RAD}{\ensuremath{\mathord{\mathcal{R}}}\xspace}

\newcommand{\FOUR}{\ensuremath{\mathord{\mathcal{F}}}\xspace}
\newcommand{\ICLOS}[1]{\ensuremath{\mathord{\overline{\overline{#1}}}}\xspace}
\newcommand{\Qn}{{\ensuremath{\mathord{\QQ^\IX{n}}}}\xspace}
\newcommand{\Dn}{{\ensuremath{\mathord{\DD^\IX{n}}}}\xspace}
\newcommand{\Cn}{{\ensuremath{\mathord{\CC^\AIX{n}}}}\xspace}
\begin{document}
\title[Asymptotic Hyperfunctions]{Asymptotic Hyperfunctions, Tempered
  Hyperfunctions, and Asymptotic Expansions}
\author[A.\ U.\ Schmidt]{Andreas U.\ Schmidt}
\date{3rd July 2004}
\address{Fraunhofer-Institute for Secure Information Technology\\
Dolivostraße 15\\
64293 Darmstadt\\
Germany}
\urladdr{http://www.math.uni-frankfurt.de/\~{}aschmidt}
\email{\href{mailto:aschmidt@math.uni-frankfurt.de}{aschmidt@math.uni-frankfurt.de}}
\thanks{This research was supported by the Deutsche
Forschungsgemeinschaft \href{http://www.dfg.de}{DFG}.
The author is very grateful for the help of and the discussions
with E.~Brüning, F.~Constantinescu, R.~Estrada and S.~Nagamachi. The hospitality
of the Department of Mathematics and Applied Mathematics, University of 
Durban--Westville, South Africa, the 
\href{http://www.df.unipi.it}{Dipartimento di Fisica E.~Fermi},
Universit\`{a} di Pisa, and the 
\href{http://www.infn.it}{INFN}, 
sezione Pisa, Italy, is gratefully acknowledged}
\begin{abstract}
  We introduce new subclasses of Fourier hyperfunctions of mixed type,
  satisfying polynomial growth conditions at infinity, and develop their
  sheaf and duality theory. We use Fourier transformation and duality to
  examine relations of these \emph{asymptotic} and \emph{tempered}
  hyperfunctions to known classes of test functions and distributions,
  especially the Gelfand--Shilov-Spaces.  Further it is shown that the
  \emph{asymptotic} hyperfunctions, which decay faster than any negative
  power, are precisely the class that allow asymptotic expansions at
  infinity. These asymptotic expansions are carried over to the
  higher-dimensional case by applying the \emph{Radon transformation}
  for hyperfunctions.
\end{abstract}
\maketitle
%
\section{Introduction}
Since the advent of Sato's hyperfunctions~\cite{SAT59,SAT60}, and
the introduction of Fourier hyperfunctions by Kawai~\cite{KAW70},
the research field of hyperfunctions has become grossly
diversified. Main branches are the algebro-analytic~\cite{b:KKK86}
and the functional analytic approach to the subject. Within the
latter, in which the present study takes its place, a large
number of special classes of hyperfunctions has been considered,
cf.\ the introductions of~\cite{ITO88,ITO92,SAB85}. The
construction of two new subclasses of Fourier hyperfunctions in
this article is driven by two motives: Firstly, their relation to
known classes of distributions and hyperfunctions, and, secondly
yet not less, their intended application. The two classes of
\emph{tempered} respectively \emph{asymptotic} hyperfunctions that
we consider, satisfy two extreme cases of polynomial bounds at
infinity. The latter fall off faster than any power, while the
former are allowed to grow as an arbitrary finite power. With
respect to the first motive above, we show that tempered and
asymptotic hyperfunctions fit into and extend the scheme of
generalized functions introduced by Gelfand and
Shilov~\cite{b:GS64}. In this way, we gain insight in the
operations of duality and Fourier transform on our and several
other spaces of test and generalized functions, paralleling
earlier studies~\cite{SEB58,HAS61,PM73}. The second motive has
two roots: The application of hyperfunctions in theoretical
physics, and the more general and classical subject of asymptotic
expansions~\cite{b:WW52,b:WAS65,b:OLV97}. To the first, there is
a long standing view that in a fundamental formulation of quantum
field theory, the mathematical problems of QFT can bee seen as a
problem of the choice of the `right' class of generalized
functions for the representation of quantum
fields~\cite{WIG81,JS88,BN89,MS92,b:STR93,SCH97A,SCH97B,SOL97A,SOL97B}.
Among other developments, this has led to formulations of QFT in
terms of ultradistributions~\cite{CT74,CT79,PIE88} and finally
hyperfunctions~\cite{NM76A,NM76B,NM77,NM79,NM86,BN89,BN98}.
Furthermore, there are results which relate these problems, and
especially the most difficult subject of
renormalization~\cite{b:HOO94}, to asymptotic
expansions~\cite{BS74,BS75,b:VDZ90}. 
Also, infrared divergences show a connection to these~\cite{SCH02}.
This altogether inspired our
interest in the possibility of asymptotic expansions for a
suitable class of hyperfunctions, and we are able to show that
our asymptotic hyperfunctions are well suited in this respect.

This article contains some of the essential parts of~\cite{b:SCH99} in
Sections~\ref{sect:struct} and~\ref{sect:asympt}, and is organized
as follows:

In Section~\ref{sect:ito} we establish the sheaf theory of
tempered and asymptotic hyperfunctions of general type with
values in a Fr{é}chet space by the duality method. The strategy
follows coarsely the proceeding of~\cite{ITO88} and uses methods
and arguments from other sources, see\eg~\cite{NAG81A,PET84}, which are
almost classical. Therefore, to omit superfluous repetitions and
in order to clarify the line of argument, we state only the core
results, postponing all proofs to Appendix~\ref{app:proofs}. A
further generalization to tempered and asymptotic hyperfunctions
with values in a general Hilbert space as in~\cite{ITO92} seems
possible, but we do not undertake this.

It should be noted at this point that in the one-dimensional case,
the sheaf theory for tempered and asymptotic hyperfunctions can
be built upon relatively elementary complex-analytic methods, as
in the case of ordinary hyperfunctions,
see~\cite{b:KAN88,b:MOR93}. 
In essence, this amounts to analogies of Runge's approximation theorem
and Mittag--Leffler's theorem with polynomial growth conditions.
This is done in~\cite{b:SCH99}, where polynomial bounds at infinity 
for hyperfunctions in one dimension are established.
For the duality theory of these hyperfunctions in one dimension,
it is useful to follow the spirit of the famous
Phragm\'{e}n--Lindelöf-Principle, to obtain polynomial bounds on
integrals along a contour around an unbounded domain in \CC from
bounds in the interior. This result of separate interest is
contained in~\cite{SCH03}. 

Section~\ref{sect:struct} explores the functional analytic structure of the
spaces of tempered and asymptotic hyperfunctions. To that end, we combine
duality to test function spaces with behavior under Fourier transformation.
We are able to show the identity of tempered hyperfunctions to the dual of
the Gelfand--Shilov-Space \SIDX{}{1}{}, see~\cite{b:GS64}. This way, we extend
the Gelfand--Shilov-Scheme of test function and distribution spaces by
hyperfunctions with polynomial growth conditions.

Section~\ref{sect:asympt} contains an application of asymptotic hyperfunctions
which we regard as essential. We use them to extend the asymptotic expansions of
distributions exhibited in~\cite{EK90,b:EK94} to hyperfunctions, 
cf.\ also the related results for ultradistributions in~\cite{CIO90}.
It turns out that
the asymptotic ones are the natural objects in the category of
hyperfunctions for such expansions. We start by exploring the one-dimensional
case. Generalization to higher dimensions could trivially be done using
cartesian
products, but we prefer a more symmetric approach which uses the Radon
transformation for hyperfunctions described by Kaneko and Takiguchi
in~\cite{KT95}.

Finally, the Appendix contains the proofs of the statements in Section~\ref{sect:ito}.
\section{Hyperfunctions by the Duality Method}
\label{sect:ito}
At the heart of the duality method lies a general theorem on the
existence and uniqueness of a flabby sheaf of Fréchet spaces under
quite weak conditions which are streamlined for the use with duals of
appropriate test function spaces. It first appeared in~\cite{JUN78} and was
further generalized in~\cite{ITO88}. We use a slightly weaker statement,
which is sufficient for our needs.
\begin{theo}[Shapira--Junker--Ito, {\cite[{Theorem 1.2.1}]{ITO88}}]
\label{theo:SJI}
  Take $X$ to be a locally compact, $\sigma$-compact  topological
  space satisfying the second axiom of countability. Assume there is a
  collection $\{F_K\}$ of Fréchet spaces labeled by the compact
  subsets $K$ of $X$ such that $F_\varnothing=0$ and for any two
  compacta $K_1$, $K_2\subset X$ holds:
\begin{enumerate}
\item If $K_1\subset K_2$, then exists a continuous injection
  $i_{K_1,K_2}\colon F_{K_1}\to F_{K_2}$.
\item If $K_1\subset K_2$ is such that every connected component of
  $K_2$ intersects $K_1$, then $i_{K_1,K_2}$ has dense image.
\item The sequence of Fréchet spaces
\[
\begin{CD}
  0 @>>> F_{K_1\cap K_2} @>>> F_{K_1}\oplus F_{K_2} @>\lambda>>
  F_{K_1\cup K_2} @>>> 0,
\end{CD}
\]
with $\lambda\colon(u_1,u_2)\mapsto u_1-u_2$, is an exact topological
sequence.
\item Fore every at most countable family $\{K_i\}$ of compacta in
$X$ holds $F_K=\bigcap_i F_{K_i}$, where $K=\bigcap_i K_i$.
\end{enumerate}
Then, there exists exactly one flabby sheaf $\mathcal{F}$ on $X$ with
$\Gamma_K(X,\mathcal{F})=F_K$ for every compact set $K$ in $X$.
\end{theo}
In practice, the spaces $F_K$ will be spaces of locally analytic
func­tio­nals on real subsets of $\CC^n$. To obtain all types
of Fourier hyperfunctions, the base space $X$ is set out as a
combination of $\RR^n$ and two types of radial compactifications
of $\RR^n$: As usual we denote by $\DD^n$, $n\in\NN$, the radial
compactification of $\RR^n$ in the sense of Kawai,
see\eg\cite{b:KAN88}.  To denote the base spaces on which germs
of holomorphic functions respectively\ hyperfunctions of
arbitrary mixed type live, we use triple indices of nonnegative
integers $\IX{n}\DEF(n_1,n_2,n_3)$, $\IX{n}\in\II$. Here, we
denote by \II the subset of $\IX{n}\in\NN_0^3$ such that
$\AIX{n}\DEF n_1+n_2+n_3\neq0$. With this, we set
$\Qn\DEF\CC^{n_1}×(\DD+i\RR)^{n_2}×(\DD× i\DD)^{n_3}$,
for $\IX{n}\in\II$. Here, the real subspace
$\Dn\DEF\RR^{n_1}×(\DD^1)^{n_2}×(\DD^1)^{n_3}$ is conceived as a
compact subset of $\Qn$. We will later introduce separate symbols
for the common special cases of indices $\IX{n}$ corresponding to
ordinary, Fourier-, modified, and mixed type hyperfunctions. The
reader will find it easy to reconstruct the notation
of~\cite[{Section 2.1}]{ITO88} from ours. We set $z=(z',z'',z''')$
for $z\in\Cn$, with $z'=(z_1,\ldots,z_{n_1})$,
$z''=(z_{n_1+1},\ldots,z_{n_1+n_2})$,
$z'''=(z_{n_1+n_2+1},\ldots,z_\AIX{n})$. For any $S\subset\Qn$ we
write $S_\Cn$ for $S\cap\Cn$. We denote by $\ICLOS{U}$ the
closure of $U$ in $\Qn$ and by $K^\circ$ the interior of $K$.
\begin{defi}
  For an open set $U\subset\Qn$ let \HOPI{U} (resp. \HOMI{U}) be the
  space of all holomorphic functions $f$ on $U_\Cn$, such that for any
  compact set $K\subset U$ there exists a $\gamma\in\RR$ (resp. for all
  $\gamma\in\RR$) and
\[
\sup_{K_\Cn}\ABS{f(z)(1+\ABS{ z''}+\ABS{z'''})^{-\gamma}}<\infty
\]
holds.
The \textbf{sheaves \HOPMI{} of germs of tempered
  \textnormal{respectively\ }asymptotic holomorphic functions} are the
sheafifications of the pre-sheaves generated by the spaces
\HOPMI{U} of local sections.
\end{defi}
Next we introduce topologies on the spaces of local sections.
\begin{defi}
\label{defi:spaces}
  Let $K\subset\Qn$ be compact and $U\subset\Qn$ be open.  For $m\in\ZZ$
  and a compact set $K\subset\QQ^\IX{n}$ we set
\begin{align*}
  \NO{f}{m}{K}\DEF \sup_{K_\Cn}\ABS{f(z)(1+\ABS{
  z''}+\ABS{z'''})^{-m}},
\end{align*}
whenever this makes sense for a function $f$.  We denote by \HOBO[m]{K}
the space of holomorphic functions $f$ on $K^\circ_\Cn$, which are continuous on
$K_\Cn$, and such that $\NO{f}{m}{K}<\infty$ holds. Choose a fundamental
system $\{V_m\}$ of neighborhoods of $K$ with $V_{m+1}\Subset V_m$ and
a sequence $\{L_m\}$ of compacta which exhausts $U$. We set
\begin{align*}
  \HOPI{K}&\DEF\varinjlim\HOBO[m]{\ICLOS{V_m}}, &
  \HOMI{U}&\DEF\varprojlim\HOBO[-m]{L_m}, \\
  \HOPI{U}&\DEF\varprojlim\HOPI{L_m}, &
  \HOMI{K}&\DEF\varinjlim\HOMI{V_m},
\end{align*}
thereby introducing locally convex topologies on these spaces.
\end{defi}
\begin{prop}
\label{prop:DFS}
  The spaces \HOPMI{K} are DFS-spaces and \HOPMI{U} are
  FS-spaces. All these spaces are nuclear.
\end{prop}
The \textbf{sheaves of germs of tempered \textnormal{respectively\ } asymptotic
  real analytic functions} are defined by
$\PPMI{}\DEF\HOPMI{}|_\Dn$. The spaces of sections
\PPMI{K} of \PPMI{} on a compact set
$K\subset\Dn$ are the DFS-spaces \HOPMI{K}.
The spaces of local sections of the sheaves \HOPMI{} and
\PPMI{} exhibit the usual tensor product decomposition property:
\begin{prop}
\label{prop:tensor}
  For compact sets $K\subset\Qn$ and $L\subset\QQ^\IX{m}$, we have the following
  topological isomorphisms:
\begin{enumerate}
\item $\HOPMI{U× V}\cong
  \HOPMI{U}\CTENS\HOPMI{V}$, $U\subset\Qn$,
  $V\subset\QQ^\IX{m}$ open.
\item $\HOPMI{K× L}\cong
  \HOPMI{K}\CTENS\HOPMI{L}$, $K\subset\Qn$,
  $L\subset\QQ^\IX{m}$ compact.
\item $\PPMI{K× L}\cong
  \PPMI{K}\CTENS\PPMI{L}$, $K\subset\Dn$,
  $L\subset\DD^\IX{m}$ compact.
\item $\PPMI{\QQ^{(n_1,n_2,n_3)}}\cong
  \mathcal{A}(\RR^{n_1})\CTENS\PPMI{\DD^{(0,n_2,n_3)}}$,
\end{enumerate}
where $\mathcal{A}(\RR^{n})$ denotes the space of ordinary real analytic
functions on $\RR^n$.
\end{prop}
By duality, one could derive Schwartz-type kernel theorems for the
tempered and asymptotic hyperfunctions to be defined below from
the above proposition, as in~\cite[Section 3.1]{ITO88}
or~\cite{CKL95}, but we will be content with leaving this issue on
the level of test functions.
\begin{theo}
\label{theo:H1vanish}
  For every compact set $K\subset\Dn$ holds
  $\RELCO{1}{K}{\Dn}{\PPMI{}}=0$.
\end{theo}
This theorem is the basis for the localization of hyperfunctions.
Name­ly, by considering the long exact sequence of cohomology groups
\begin{align*}
    0 &\longrightarrow \PPMI{K_1\cup K_2}
      \longrightarrow \PPMI{K_1}\oplus\PPMI{K_2}
      \longrightarrow \PPMI{K_1\cap K_2}\\
      &\longrightarrow \RELCO{1}{K_1\cup K_2}{\Dn}{\PPMI{}}
      \longrightarrow \cdots
\end{align*}
for two compact sets $K_1$, $K_2\subset\Dn$,
we immediately derive from it the following important conclusion, which
is dual to condition iii) of Theorem~\ref{theo:SJI}:
\begin{coro}
\label{coro:locseq}
  The following sequence is exact:
\[
    0 \to \PPMI{K_1\cup K_2}
      \longrightarrow \PPMI{K_1}\oplus\PPMI{K_2}
      \longrightarrow \PPMI{K_1\cap K_2}\to 0.
\]
\end{coro}
The last main ingredient is an approximation theorem of Runge type.
\begin{theo}
\label{theo:runge}
 \PPMI{\Dn} is dense in \PPMI{K} for $K\subset\Dn$
 compact.
\end{theo}
Now let $E$ be any Fréchet space. For an open set $U\subset\Qn$
we call the spaces $\HOPMID{U;E}\DEF L(\HOPMI{U};E)$ of all
continuous linear mappings from \HOPMI{U} into $E$ the
\textbf{asymptotic} respectively \textbf{tempered analytic
functionals on $U$ with values in $E$}. Similarly, we define
$\HOPMID{K;E}\DEF L(\HOPMI{K};E)$ respectively $\PPMID{K;E}\DEF
L(\HOPMI{K};E)$ for $K\subset\Qn$ respectively $K\subset\Dn$
compact. All these spaces are endowed with the topology of
convergence on compact subsets. Then, by virtue
of~\cite[Proposition 50.5]{b:TRE67}, see also~\cite[proof of
Theorem 5.7]{NAG81B}, we have:
\begin{prop}
  For any Fréchet space $E$ holds:
  \begin{enumerate}
  \item $\HOPMID{U;E}\cong\HOPMID{U}\CTENS E$, for $U\subset\Qn$ open.
  \item $\HOPMID{K;E}\cong\HOPMID{K}\CTENS E$, for $K\subset\Qn$ compact.
  \item $\PPMID{K;E}\cong\PPMID{K}\CTENS E$, for $K\subset\Dn$ compact.
  \end{enumerate}
\end{prop}
We say that a compact set $K\subset U\subset\Qn$ is a \textbf{carrier}
for a section $F\in\HOPMID{U;E}$ if $F$ can be extended to an element of
\HOPMID{K;E}. The functional $F$ is said to be carried by an open subset $V$
in $U$ if it is carried by some compact subset of $V$. If a compact set
$K$ in $U\subset\Dn$ has the Runge property, and thus \PPMI{U} is dense in
\PPMI{K} by Theorem~\ref{theo:runge}, then $F$ is carried by
$K$ if and only if it is carried by all open neighborhoods of $K$ in
$U$.

By using the dual of the exact sequence of Corollary~\ref{coro:locseq}, and
by induction, we easily see that $\bigcap_i\PPMID{K_i;E}=\PPMID{\bigcap_i
  K_i;E}$ for every countable family $\{K_i\}$ of compacta in
$\Dn$. Then, Zorn's Lemma implies that for every functional
$F\in\PPMID{\Dn;E}$ with $F\neq0$, we can find a smallest compact
set $K$ in \Dn which is a carrier for $F$. We call $K$ the
\textbf{support} of $F$ and denote it by $\supp(F)$,
cf.~\cite[Theorems 2.3.4 and 2.3.5]{ITO88}. Then, the identity
$\PPMID{K_i;E}=\{F\in\PPMID{\Dn;E}\mid\supp(F)\subset K\}$ easily
follows.

With these preparations, we are ready to define the sheaves of
\textbf{tempered} and \textbf{asymptotic hyperfunctions of
general type with values in a Fréchet space $E$}. Namely, the
mapping $K\mapsto\PPMID{K;E}$, which assigns a Fréchet space to
every compact set $K\subset\Dn$, satisfies all conditions of the
Shapira--Junker--Ito-Theorem \ref{theo:SJI}, cf.~\cite[proof of
Theorem 2.4.1]{ITO88}. Thus we have:
\begin{theo}
\label{theo:sheaf}
  There exists exactly one flabby sheaf ${}^E\!\QPMI{}$ such that for
  every compact set $K\subset\Dn$ holds
  $\Gamma_K(\Dn,{}^E\!\QPMI{})=\PIDX{\mp\infty}{\prime}{K;E}$.
\end{theo}
There is a natural embedding of flabby sheaves
\[
{}^E\!\SHQMI\hookrightarrow{}^E\!\SHQPI\hookrightarrow{}^E\!\SHQ[]
\]
of asymptotic into tempered into ordinary Fourier hyperfunctions on \Dn,  induced by the continuous inclusions of the respective test functions spaces.
\section{The Structure of Tempered and Asymptotic Hyperfunctions}
\label{sect:struct}
In this section, we specialize to the case of scalar-valued,
unmodified, tempered and asymptotic hyperfunctions\ie we consider
$\QPMI{}={}^\CC\!\QPMI{}$ on
$\DD^n=\DD^{(0,n,0)}\subset\QQ^n=\QQ^{(0,n,0)}$. The following
theorem establishes the orthantic boundary value representation
for global sections of these sheaves, cf.~\cite[§7.1]{b:KAN88}.
We do not go into developing a duality theory for local sections,
resembling Poincarè--Serre duality for cohomology groups, but
rather present duality theorem which relates globally defined
tempered and asymptotic hyperfunctions to boundary values of
holomorphic functions with the same growth or decay behavior.
\begin{theo}
There is a linear, topological isomorphism
\[
\QPMI{\DD^n}\cong
\HOIDX{±\infty}{}{W\#\DD^n}\Bigm/\sum_{j=1}^n
\HOIDX{±\infty}{}{W\#_j\DD^n},
\]
for every open, cylindrical neighborhood $W$ of $\DD^n$ in $\QQ^n$.
\end{theo}
Here, for $W=W_1×\ldots× W_n$ and a compact, cylindrical subset
$K=K_1×\ldots× K_n$ of $W$ we define
\begin{align*}
W\#K&\DEF (W_1\setminus K_1)×\cdots×(W_n\setminus K_n)\\
\intertext{and}
W\#_jK&\DEF (W_1\setminus K_1)×\cdots
 ×\underbrace{W_j}_{\text{omitted}}×\cdots×(W_n\setminus K_n).
\end{align*}
We will not give a detailed proof of this theorem, since the
existing ones for
Fourier hyperfunctions can be literally applied in our case, 
see\eg the clear exposition in~\cite[Part~C]{BN89}. 
Let us nevertheless comment on the essential points.
For an equivalence class $[F]$ in one of the quotients defined
above and a function $f\in\PIDX{\mp\infty}{}{\DD^n}$ one defines
an inner product
\[
\SPROD{[F]}{f}=
-\int\limits_{\Gamma_1}\cdots\int\limits_{\Gamma_n}
F(z_1,\ldots,z_n) f(z_1,\ldots,z_n) \dd z_1\cdots \dd z_n,
\]
where the integration plane $\Gamma=\Gamma_1×\ldots×\Gamma_n$ has to be chosen to lie in the common domain of holomorphy of $F$ and $f$. Since it is of the form of a cartesian product, Cauchy's Theorem ensures independence of the bilinear form of the special choice of $\Gamma$.
One easily sees that the linear functional $T_{[F]}=\SPROD{[F]}{\cdot}$ is continuous. That the mapping $F\mapsto T_{[F]}$ is injective is essentially an application of Cauchy's integral formula, but with an exponentially decaying kernel. This kernel
\[
h_z(w)\DEF\prod_{i=1}^{n}\MTPI\cdot\frac{\ee^{-(z_i-w_i)^2}}{z_i-w_i},
\]
is also used to show surjectivity of $T$ by evaluating it on a given functional 
$T\in\SQ[\PMI]{\DD^n}$. The function $T(h_z)$ is in \HOIDX{±\infty}{}{W\#\DD^n}, 
since $h$ preserves the present asymptotic resp.\ tempered growth condition as can
easily be verified by explicit estimation, see~\cite{b:SCH99}, 
and defines $T$ via $T_{[T(h_z)]}=T$ 
and accordingly, $T(h_z)$ is called a \textbf{defining function} for the 
hyperfunction $T$.

The set $W\#\DD^n$ decomposes into $2^n$ connected components labeled by the signs $\sigma=(\sigma_1,\ldots,\sigma_n)$ of the imaginary parts of the components $(z_1,\ldots,z_n)$ of the coordinate $z$. By this decomposition, every tempered or asymptotic hyperfunction possesses the \textbf{orthantic boundary value representation}
\[
f(x)=\sum_\sigma F_\sigma(x+i\Gamma_\sigma 0),\
F_\sigma(z)\in\HOPMI{\DD^n+i\Gamma_\sigma 0},
\]
where $\Gamma_\sigma\DEF\{x\in\RR^n\mid \sigma\cdot x> 0\}$ is
the \textbf{$\sigma$-th orthant} and $F_\sigma$ is holomorphic on
an \textbf{infinitesimal wedge of type $\DD^n+i\Gamma_\sigma 0$},
see~\cite[Page 82]{b:KAN88}. Contact with representations by
boundary values from other infinitesimal wedges can be made by
convolution of a hyperfunction $f\in\QPMI{\DD^n}$ with the
\textbf{exponentially decreasing Radon decomposition kernel}
$W_\ast(x,\omega)$, $\omega\in S^{n-1}$, which preserves the
polynomial growth conditions on $f$, see~\cite[Appendix]{KT95}.
Consistency of all such representations is assured by
\textbf{Martineau's edge of the wedge theorem with polynomial
decay conditions}, which we may cite now in a form suitable for
tempered and asymptotic hyperfunctions:
\begin{theo*}[{\cite[Theorem A.6]{KT95}}]
  Let $f(x)$ be a Fourier hyperfunction with a set of defining functions
  $\{F_j(z)\in\HOPMI{\DD^n+i\Gamma_j0}\}_{j=1}^N$. Assume $f=0$ in \QPMI{\DD^n}. Then for
  any choice of proper sub-cones $\Gamma'_j\Subset\Gamma_j$ there exist 
 wedge-analytic functions
  $F_{jk}\in\HOPMI{\DD^n+i(\Gamma'_j+\Gamma'_k)0}$ such that
  \begin{displaymath}
    F_{jk}=- F_{kj},\ F_j(z)=\sum_{k=1}^N F_{ik}(z)
  \end{displaymath}
on an infinitesimal wedge of type $\DD^n+i\Gamma_j0$.
\end{theo*}
The Fourier transformation $\FOUR$ on \QPMI{\DD^n},
see~\cite[Chapter 7]{b:KAN88}, can be defined as usual by taking
the boundary value of the Fourier--Laplace-Transfor­ma­tion of
a single boundary value and extending linearly to the formal sums
representing $f\in\QPMI{\DD^n}$. It is consistent with the
embedding $\QPMI{\DD^n}\hookrightarrow\SQ[]{\DD^n}$ (since the
boundary value representations are), and the Fourier inversion
formula holds.

We first consider the Fourier transform of the space \SPMI{\DD^n} of asymptotic real analytic functions,
which is the test function space of \SQPI{\DD^n}. As one would expect, $\FOUR\SPMI{\DD^n}$ is a space of exponentially decreasing $C^\infty$-functions. Closer examination shows that it is one of the spaces \SIDX{\alpha}{\beta}{} introduced by Gelfand and Shilov. We give an equivalent definition of it:
\begin{defi}[{\cite[Chapter IV,\ §3]{b:GS64}}]
We set
\[
\SSS[1]{\RR^n}\DEF\Bigl\{ f\in C^\infty(\RR^n) \Bigm|
    \exists \delta>0, \forall \alpha: \sup_x\ABS{D^\alpha_x f(x)}\ee^{\delta\ABS{x}}
    <\infty \Bigr\}.
\]
The topology of  \SSS[1]{\RR^n} is that of an inductive limit
\[
\SSS[1]{\RR^n}=\mathop{\varinjlim}_{A\to\infty}
\SIDX{1,A}{}{\RR^n},
\]
of countably normed spaces, where
$\SIDX{1,A}{}{\RR^n}$ is the space of all infinitely differentiable functions $f$ for which all the norms
\[
  \NORM{f}_{\SSS[1]{};m,p}\DEF
  \sup_{x,\ABS{\alpha}\leq p} \ABS{D^\alpha_x
    f(x)}\ee^{m^{-1}(1-p^{-1})\ABS{x}},
\]
are finite, where $m=\ee A$, and $p=2,3,\ldots$.
\end{defi}
\begin{theo}
\label{theo:PW1}
The Fourier transformation
$\FOUR\colon \SPMI{\DD^n}\lirarrow\SSS[1]{\RR^n}$ induces a linear topological isomorphism.
\end{theo}
\begin{proof}
We choose an equivalent representation of the space \SPMI{\DD^n} as an inductive limit of countably normed spaces:
\begin{displaymath}
  \SPMI{\DD^n}\cong\mathop{\varinjlim}\limits_m\HOBO[-\infty]{\ICLOS{U_m}}\text{, with }
  \HOBO[-\infty]{\ICLOS{U_m}}\DEF\mathop{\varprojlim}\limits_k\HOBO[-k]{\ICLOS{U_m}}.
\end{displaymath}
Here we use the special system $U_m=\DD^n+i\{\ABS{\IM z}<1/m\}$ of neighborhoods of $\DD^n$. Now let $f\in\HOBO[-\infty]{\ICLOS{U_m}}$ for an $m\in\NN$.
Since $f$ is an asymptotic function on the whole domain $\ICLOS{U_m}_{\CC^n}$ we can use Cauchy's Theorem and dominated convergence to calculate its Fourier transform $\widehat{f}$ by shifting the integration plane as follows:
\[
D^\alpha_\xi\widehat{f}(\xi)=\int\limits_{z=\RR^n+iy}(-iz)^\alpha \ee^{-iz\xi}f(z) \dd z,
\]
with arbitrary $\ABS{y}\leq1/m$. Choosing $z=x\mp i/m$ for $\xi\gtrless 0$, we estimate
 \begin{align*}
   \ABS{\smash{D^\alpha_\xi\widehat{f}(\xi)}}&\leq
    \NORM{f}_{-\ABS{\alpha}-n-1,\ICLOS{U_m}}\cdot
      \ABS{\int_{\RR^n} \ee^{-ix\xi}\ee^{-\ABS{\xi}/m}(1+\ABS{x})^{-n-1} \dd x}\\&\leq
    C_n\cdot\NORM{f}_{-\ABS{\alpha}-n-1,\ICLOS{U_m}}\cdot \ee^{-\ABS{\xi}/m},
\end{align*}
with certain $C_n>0$. This shows that $\widehat{f}$ is an exponentially decreasing $C^\infty$-function\ie an element of \SSS[1]{\RR^n}. Thus, for the norm \NODX{\cdot}{\SSS[1]{};m}{p} on \SSS[1,A]{\RR^n} with $A=m/\ee$ we have
\begin{align*}
  \NODX{\smash[t]{\widehat{f}}}{\SSS[1]{};m}{p}&\leq
  C\cdot \NORM{f}_{-p-n-1,\ICLOS{U_m}}\cdot \sup_{\xi}
  \ee^{-\ABS{\xi}/m} \ee^{(1-p^{-1})\ABS{\xi}/m} \\
  &\leq C\cdot \NORM{f}_{-p-n-1,\ICLOS{U_m}}\cdot \sup_{\xi}
  \ee^{-\ABS{\xi}/(mp)} \leq C\cdot \NORM{f}_{-p-n-1,\ICLOS{U_m}},
\end{align*}
and similarly for every $m'>m$, which shows continuity of $\FOUR$ with respect to the inductive limit topologies of \SPMI{\DD^n} and \SSS[1]{\RR^n}. Then, by the classical Fourier inversion formula, $\FOUR$ is a continuous linear bijection and the continuity of $\FOUR^{-1}$ follows similarly as above.
\end{proof}
Since the Fourier transformation acts on the Gelfand--Shilov-Spaces by exchanging the indices, we can immediately place \SPMI{\DD^n} itself into the \SIDX{\alpha}{\beta}{}-scheme. Tempered hyperfunctions also find their place, since \SPMI{\DD^n}  are their test functions:
\begin{coro}
\label{coro:PW1}
The space \SPMI{\DD^n} is topologically isomorphic to the function space \SIDX{}{1}{\RR^n} and \SQPI{\DD^n} is topologically isomorphic to the space \SIDX{}{1\prime}{\RR^n}. The Fourier transformation induces a mapping
$\FOUR\colon \SQPI{\DD^n}\lirarrow\SIDX{1}{\prime}{\RR^n}$, which is a
linear, topological isomorphism.
\end{coro}
Modeled after the scheme exhibited above, we can now examine the Fourier transform of asymptotic hyperfunctions. It is by now clear that $\FOUR\SQMI{\DD^n}$ is a space of smooth functions. They exhibit the infra-ex­po­nen­ti­al growth property which is typical for Fourier transforms of Fourier hyperfunctions.
\begin{defi}
\label{defi:CIIE}
We denote by  \CIIE{\RR^n} the space of
\textbf{infra-exponential},  smooth functions. These are all $f\in
C^\infty(\RR^n)$ such that: For all $k\in\NN$ and $\EPS>0$
exists a constant $C_{k,\EPS}>0$ with
\[
\ABS{\frac{\partial^\alpha f(x)}{\partial x^\alpha}}\leq
C_{k,\EPS}\ee^{\EPS\ABS{\RE z}},
\]
for all multi-indices $\alpha\in\NN^n$ with $\ABS{\alpha}\leq k$. We equip \CIIE{\RR^n} with the topology of a countably normed space induced by the norms
\begin{displaymath}
  \NODX{f}{\CIIE{};m}{p}\DEF\sup_{x,\ABS{\alpha}\leq p} \ABS{D^\alpha_x
    f(x)}\ee^{-\ABS{x}/m}
\end{displaymath}
for all $p$, $m\in\NN$, and consider \CIIE{\RR^n} to be completed in this topology.
\end{defi}
We note aside that with this topology, \CIIE{} is isomorphic to the test function space which was denoted by $\mathcal{P}$ in \cite{EK90}.
\begin{theo}
\label{theo:PW2}
The Fourier transformation
$\FOUR\colon \SQMI{\DD^n}\lirarrow\CIIE{\RR^n}$ induces a linear topological isomorphism.
\end{theo}
\begin{proof}
Since \FOUR extends to linearly to sums of boundary values, it suffices to consider an asymptotic hyperfunction $f$ represented by a single boundary value $f(x)=F(x+i\Gamma_\sigma0)$ from an orthant $\Gamma_\sigma$. So, let $W$ be an infinitesimal wedge of type $\DD^n+i\Gamma_\sigma0$ and $F(x+iy)\in\HOMI{W}$ be a defining function for $f$. Then, for every compact set $L\subset\RR^n$ such that $K=\DD^n+i L$ is compact in $W$, holds the estimate:
\begin{equation}
\begin{aligned}
  \ABS{\smash{D^\alpha_\xi\widehat{f}(\xi)}}&=
  \ABS{\int_{\IM z=y} (-iz)^\alpha \ee^{-iz\xi}F(z) \dd z}\\
  &\leq \NODX{F}{-\ABS{\alpha}-n-1}{K}
    \ABS{\int_{\RR^n} \ee^{-ix\xi}\ee^{y\xi}(1+\ABS{x})^{-n-1} \dd x},
\end{aligned}
\tag{\dag}
\label{eq:PW2est}
\end{equation}
with arbitrary $y\in L$. This shows $\FOUR f\in\CIIE{\RR^n}$
since $\ABS{y}$ can be made arbitrarily small, leaving the
integral unchanged by Cauchy's Theorem. On the other hand, let
$\widehat{G}(\xi)\in\CIIE{\RR^n}$ be such that all derivatives of
$\widehat{G}$ decrease exponentially outside the closed cone
$\overline{\Gamma_\sigma}$, which can be achieved by eventually
decomposing the original function utilizing exponentially
decreasing multipliers. Then, the inverse Fourier transform
$G=\FOUR^{-1}\widehat{G}$ is a boundary value
$G(z)\in\HOIE{\DD^n+i\Gamma_\sigma0}$ and thus a Fourier
hyperfunction, cf.~\cite[Proposition 8.3.2]{b:KAN88}. It is an
easy calculation to show $G(z)=O(\ABS{\RE z}^{-\infty})$ locally
uniformly in $\IM z$. This shows $\FOUR^{-1}
\widehat{G}\in\SQMI{\DD^n}$. The inversion formula for Fourier
hyperfunctions, see~\cite[Theorem 8.3.4]{b:KAN88}, thus implies
that $\FOUR$ is a linear bijection from \SQMI{\DD^n} onto
\CIIE{\RR^n} with inverse $\FOUR^{-1}$. It remains to show
continuity. If we choose a special exhausting sequence of compacta
$\{K_j=\DD^n+i L_j\}_{j\in\NN}$ for $W$, such that some points of
the cylindrical surface $\{\ABS{\IM z}= \ABS{y}= 1/j\}$ are
contained in $K_{j\CC^n}$, then we can make $y$ in
estimate~\eqref{eq:PW2est} small enough to conclude
\begin{displaymath}
  \ABS{\smash{D^\alpha_\xi\widehat{f}(\xi)}}\leq
  C\cdot\NODX{F}{-\ABS{\alpha}-n-1}{K_j}\cdot \ee^{\ABS{\xi}/j}.
\end{displaymath}
This yields
\begin{displaymath}
    \NODX{\smash[t]{\widehat{f}}}{\CIIE{};m}{p} \leq
    C\cdot \NODX{F}{-p-n-1}{K_j} \sup_\xi
    \ee^{(j^{-1}-m^{-1})\ABS{\xi}}\leq
    C\cdot \NODX{F}{-p-n-1}{K_j},
\end{displaymath}
for all $j\geq m$,
showing continuity of \FOUR in the topologies of \SQMI{\DD^n} and \CIIE{\RR^n}.
\end{proof}
Again, we can immediately draw the following conclusion:
\begin{coro}
The Fourier transformation
$\FOUR\colon \SPPI{\DD^n}\lirarrow\CIIED{\RR^n}$ induces a linear topological isomorphism.
\end{coro}
The test function spaces \SIDX{\alpha}{\beta}{} for $\alpha$,
$\beta\geq1$, are ordered in the Gelfand--Shilov scheme according
to two characteristics: Growth order, controlled by the lower
index $\alpha$, which ranges from exponential decay for
$\alpha=1$ to rapid (asymptotic in our terminology) decay for
$\alpha=\infty$. And regularity, which is that of real analytic
functions in strip-like neighborhoods of the real axis for
$\beta=1$\ie the typical regularity of test function spaces of
Fourier hyperfunctions, and on the other hand simple
$C^\infty$-functions for $\beta=\infty$ (of course satisfying the
growth conditions demanded by $\alpha$ to all derivatives). Note
that\eg the Schwartz-Space \SSS{} is nothing but
\SIDX{\infty}{\infty}{}, and it is long known,
see~\cite[Proposition 2.1]{NM76A}, that $\SIDX{1}{1}{\RR^n}$ is
exactly the space $\PIDX{\ast}{}{\DD^n}$ of exponentially
decaying, real analytic test functions whose dual is the space
\SQ[]{\DD^n} of Fourier hyperfunctions.

We can use our two Paley--Wiener-Type Theorems~\ref{theo:PW1} and~\ref{theo:PW2} and their corollaries to extend this scheme largely to include asymptotic, tempered and Fourier hyperfunctions.  The result is shown in the figure
above.
On the ordinate are marked four types of growth conditions: Exponential decay, asymptotic\ie rapid decay, tempered growth and infra-exponential growth, symbolized in that order by \HOIDX{\ast}{}{}, \SHOMI, \SHOPI, \HOIE{}. Here, the use of the symbol for holomorphic functions is justified by the embedding of the various spaces into spaces of hyperfunctions with the namely growth conditions\ie the existence of representations by boundary values of holomorphic functions exhibiting these conditions. The regularities marked on the abscissa are \SHP[] for functions which are real analytic in strip-neighborhoods, $C^\infty$ for smooth functions, \SDDFO[]{} for distributions of finite order, and \SHQ[]{} for hyperfunctions.
\begin{figure}[t!]
\begin{small}
\begin{displaymath}
\psset{xunit=1.9cm}
\psset{yunit=1.9cm}
\begin{pspicture}(-1.5,-1.5)(3.5,3.5)
\uput{2mm}[180](-1.5,1.15){\rotateleft{Singularity}}
\uput{2mm}[270](1,-1.5){Growth}
\psaxes[labels=none,ticks=none,linewidth=1.2pt]{->}(-1.5,-1.5)(3.5,3.5)
\multips(0,-1.5)(1,0){4}{\psline[linewidth=1.2pt](0,0)(0,0.1)}
\multips(-1.5,0)(0,1){4}{\psline[linewidth=1.2pt](0,0)(0.1,0)}
\psframe*[linecolor=mid1](0,0)(1,1)
\psline[linestyle=dotted](-1.5,-1.5)(-0.1,-0.1)
\psline[linestyle=dotted](0.1,0.1)(0.9,0.9)
\psline[linestyle=dotted](1.1,1.1)(3.2,3.2)
  \uput{3mm}[90](0,-1.5){\HOIDX{\ast}{}{}}
  \uput{3mm}[90](1,-1.5){\SHOMI}
  \uput{3mm}[90](2,-1.5){\SHOPI}
  \uput{3mm}[90](3,-1.5){\HOIE{}}
  \uput{3mm}[0](-1.5,0) {\SHP[]}
  \uput{3mm}[0](-1.5,1){$C^\infty$}
  \uput{3mm}[0](-1.5,2){\SDDFO[]{}}
  \uput{3mm}[0](-1.5,3){\SHQ[]}
  \rput(0,0){{$\SIDX{1}{1}{}\!\cong\!\SHP[\ast]$}}
  \rput(1,0){{$\SIDX{}{1}{}\!\cong\!\SHPMI$}}
  \rput(2,0){\SHPPI}
  \rput(3,0){$\bigl(\SHP[]\bigr)$}
  \rput(0,1){{\SSS[1]{}}}
  \rput(1,1){{\SSS{}}}
  \rput(3,1){\CIIE{}}
  \rput(0,2){\CIIED{}}
  \rput(2,2){\psframebox*{\SSD{}}}
  \rput(3,2){\SIDX{1}{\prime}{}}
  \rput(0,3){$\bigl(\SHQ[\ast]\bigr)$}
  \rput(1,3){\SHQMI}
  \rput(2,3){\SHQPI}
  \rput(3,3){\psframebox*{\SHQ[]}}
  \rput(1.5,1.5){\pscirclebox[boxsep=false,linestyle=dotted,fillstyle=solid]{$L^2$}}
  \psline{<->}(-0.75,-1.25)(-1.25,-0.75)
\rput(-1,-1){\psframebox*{\FOUR}}
\end{pspicture}
\end{displaymath}
\caption{A diagram of generalized functions.}
\label{fig:spaces}
\end{small}
\end{figure}
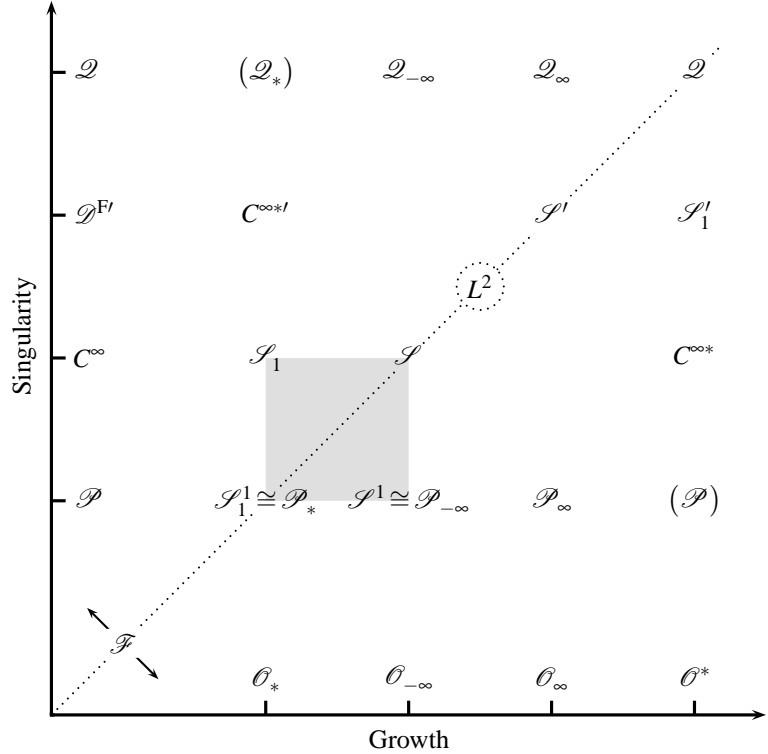

In the lower left corner of the figure we find the part of the
Gelfand--Shilov-Scheme which has been described above. The
Fourier transformation is a symmetry of the diagram which
operates by reflection on the diagonal\ie exchanging growth
conditions with singularity. The remarkable fact about the diagram
is that it incorporates a second symmetry operating by point
reflection on the center, namely duality. The combination of these two
symmetries allowed us to draw the cross-conclusions of the above
corollaries from the corresponding Paley--Wiener-Theorems. The
self-dual $L^2$ in the middle is closed under Fourier
transformation. It forms Gelfand-Triplets together with pairs of
other spaces\eg $(\SPMI{\DD^n},L^2(\RR^n),\SQPI{\DD^n})$.

We note that similar configurations of generalized functions,
which would further enhance our figure, have already been
considered by Sebasti{ã}o~e~Silva in~\cite{SEB58}, Hasumi
in~\cite{HAS61} and Park and Morimoto in~\cite{PM73}. They
considered the so-called \textbf{Fourier ultra-hyperfunctions}
$\mathcal{U}$ which are the elements of the dual space of the
space $\mathfrak{H}$ of entire functions of rapid decay. Via
Fourier transformation, $\mathfrak{H}$ corresponds to a space $H$
of smooth functions which decay faster than
$\ee^{-\gamma\ABS{x}}$ for every $\gamma>0$, and $\mathcal{U}$
corresponds to the space $\Lambda^\infty$ of
\textbf{distributions of exponential growth}.
\begin{note}
Our definition of \emph{tempered} hyperfunctions \SQPI{\DD^n}
contains an inherent ambiguity: One has to make the choice
whether the boundary value $F(x+i\Gamma0)\in\SQPI{\DD^n}$ will
have fixed growth order, say $O(x^N)$, as $y\in\Gamma$ tends to
zero or if this growth order may vary. With our definitions, the
latter is the case, for it is immediate from
Definition~\ref{defi:spaces} that $F(x+iy)$ is of a fixed growth
order in $x$ only \emph{locally uniformly} in $y\in W$, $W$ an
infinitesimal wedge of type $\Gamma$. This behavior conforms
with that of boundary value representations of tempered
distributions in $\SSD{}$. As one would expect, we find:
\begin{rem}
There is a continuous embedding
$\SPMI{\DD^n}\hookrightarrow\SSS{\RR^n}$ of test function spaces,
as
a consequence of Cauchy's estimates.
\end{rem}
The above mentioned ambiguity also appears in the case of
\emph{exponentially decreasing hyperfunctions}, cf.~\cite[Note
8.3, p.\ 411]{b:KAN88}, which explains the parentheses around
\SHQ[\ast] and \SHP[] in the figure above: In duality to the
space of infra-exponential analytic functions \SHP[] is the
relative cohomology group
\RELCO{n}{\DD^n}{\QQ^n}{\HOIDX{\ast}{}{}} which consists of
boundary values with an exponential decay at infinity that may
vary with $\IM(z)$. This was the original definition of the space
\SQ[\ast]{\DD^n} of exponentially decreasing hyperfunctions. In
contrast, to obtain the Fourier transform of \SHP[] Kaneko
defined in~\cite{KT95} exponentially decreasing hyperfunctions as
$\SQ[\ast]{\DD^n}=
\bigcup_{\EPS>0}\ee^{-\EPS\sqrt{1+x^2}}\SQ[]{\DD^n}$, which
consists of boundary values of \emph{constant} exponential decay
in $\IM(z)$. Note also that this ambiguity does not appear for
asymptotic hyperfunctions, since obviously their defining
functions remain asymptotic when approaching the real axis.
\end{note}
There is a natural relation between tempered hyperfunctions
and tempered distributions: The embedding of the corresponding
test function spaces $\SPMI{\DD^n}\hookrightarrow\SSS{\RR^n}$
remarked above has dense image~\cite[Theorem 15.5]{b:TRE67} , and
thus by duality yields the following result:
\begin{rem}
The space \SIDX{}{\prime}{\RR^n} is continuously embedded into
\SQPI{\DD^n}.
\end{rem}
The question comes up naturally, which space of distributions is
in the equivalent relation to asymptotic hyperfunctions. The
distributions \SKD{} were introduced in~\cite{GLS68} and used by
Estrada, Kanwal \textit{et~al.}\ in~\cite{EGV98} for
distributional asymptotic expansions. It will turn out that they
are related to asymptotic hyperfunctions of modified type. We
define \SKD{} and recall some of its properties from~\cite{EGV98}
and~\cite{b:EK94}.
\begin{defirem}
\label{defi:SKD}
Denote for $\gamma\in\RR$ by  \SK{\RR} the space of all
$\phi\in C^\infty(\RR)$, for which $\phi^{(k)}=O(\ABS{x}^{\gamma-k})$
holds for $\ABS{x}\to\infty$ and every $k\in\NN$.
\SK{} becomes a locally convex space equipped with the system
\begin{displaymath}
\NORM{\phi}_{K,k,\gamma}\DEF\sup_{x\in
  K}\ABS{(1+x)^{k-\gamma}\phi^{(k)}(x)}
\end{displaymath}
of seminorms for $K\subset\RR$ compact. Then
$\SK{\RR}\hookrightarrow\SK[\gamma']{\RR}$ for $\gamma\leq\gamma'$,
and we define
$
\SK[]{\RR}\smash\DEF\mathop{\varinjlim}_{\gamma\to\infty}\SK{\RR}.
$
The function algebra \SK[]{} is \textbf{normal}\ie\SSS{} is dense in \SK[]{}, and
\SK[]{\RR} is a nuclear space. We set $\SKD{\RR^n}\smash{\DEF\SKD{\RR}^{\CTENS n}}$.
\end{defirem}
To this point, we have treated
normal and modified Fourier hyperfunctions in a unified manner as far as regards
notation. From now, we use the shorthand $\SMODQMI{\DD^n}\DEF\SQMI{\DD^{(0,0,n)}}$
and $\SMODPPI{\DD^n}\DEF\SPPI{\DD^{(0,0,n)}}$ for the modified type.
\begin{prop}
There is a continuous embedding
$\SMODPPI{\DD^n}\hookrightarrow\SK[]{\RR^n}$ which induces the
embedding $\SKD{\RR^n}\hookrightarrow\SMODQMI{\DD^n}$.
\end{prop}
\begin{proof}
We need only consider the case $n=1$ due to the tensor product decomposition
property for \SMODQMI{\DD^n}, see Proposition~\ref{prop:tensor}, and the
corresponding property of \SKD{\RR^n}. If $f\in\SMODPPI{\DD}$ then
$f\in\HOBO[m]{\ICLOS{V}}$ for some $m\in\NN$ and a neighborhood
$V\subset\QQ^{(0,0,1)}$ of $\DD^{(0,0,1)}$. Without loss of generality, we can
assume that $V_\CC=\{z\in \CC\mid \ABS{\IM{z}}<\EPS(1+\ABS{\RE{z}})\}$.
 Under this conditions, it follows
from elementary properties of holomorphic functions on
wedge-shaped domains, see~\cite[Chapter 1,\ Theorem 4.2]{b:OLV97},
that $f^{(k)}\in\HOBO[m-k]{\ICLOS{U}}$ for every neighborhood
$U$ of $\DD^{(0,0,1)}$ which is relatively compact in $V$. This
shows $f\in\SK{\DD}$ for every $\gamma>m$ and thus
$\SMODPPI{\DD^n}\subset\SK[]{\RR}$. The inclusion is continuous,
since the topology of \SK[]{\RR} is weaker than that of
\SMODPPI{\DD}, which in turn is a consequence of Cauchy's
estimates. It is also finer by definition, and the inclusion has
dense image.
\end{proof}
The counterexample $z^m\ee^{iz}$ shows that the conclusion of the proposition
does not hold for asymptotic hyperfunctions of ordinary type.

To conclude this section, we prove a 
\emph{structural theorem}, as it is traditionally called in
the theory of generalized functions,  for the asymptotic and tempered hyperfunctions\ie 
that every asymptotic (tempered) hyperfunction can be rendered by applying
a certain differential operator to a continuous function with
the same asymptotic decay (tempered growth). Of course, a generic hyperfunction cannot
be yielded through a differential operator of finite order. The right
notion in this case is that of a \textbf{local (pseudo-)differential operator}.
By this we mean an infinite order differential operator $J(D)$ with
constant coefficients
\[
J(D)=\sum_{\alpha_i\geq0}b_\alpha D_x^\alpha\quad\text{with }
D_x^\alpha\DEF
\frac{\partial^{\alpha_1}}{\partial x_1^{\alpha_1}}\cdots
\frac{\partial^{\alpha_n}}{\partial x_n^{\alpha_n}},
\]
such that the coefficients satisfy the estimate
\[
\mathop{\varlimsup}_{\ABS{\alpha}\to\infty}
\sqrt[\uproot{2}\ABS{\alpha}]{\ABS{b_\alpha}\alpha!}=0.
\]
Those operators get their name from their property of preserving the supports
of hyperfunctions, thus inducing endomorphisms of the sheaf \SHQ[], and as well
of \QPMI{}. More precisely we have by combining 
Theorems~\ref{theo:sheaf},~\ref{theo:PW2}, and Corollary~\ref{coro:PW1}
with~\cite[Proposition~8.4.8]{b:KAN88}:
\begin{coro}
\label{coro:fourlocop}
  Every local operator $J(D)$ with constant coefficients induces
an endomorphism of sheafs $\QPMI{\RR^n}\longrightarrow\QPMI{\RR^n}$\ie
$J$ preserves supports, and the formula
\[
\FOUR(J(D)f)(\xi)=J(\xi)\cdot(\FOUR f)(\xi)
\]
holds for $f\in\QPMI{\RR^n}$, where in particular $J(\zeta)$ is an 
infra-exponential, entire function.
\end{coro}
\begin{theo}\labelT{struct}
  Every $f\in\QPMI{\RR^n}$ can be represented as $J(D)f_0(x)$, where
$J(D)$ is a local differential operator and $f_0\in C(\RR^n)$ is of
asymptotic decay\ie $f(x)=O(\ABS{x}^{-\infty})$ if $f\in\SQMI{\RR^n}$ and
$f_0$ is of tempered growth\ie $f(x)=O(\ABS{x}^{r})$, for some $r\in\RR$, 
if $f\in\SQPI{\RR^n}$.
\end{theo}
The proof goes coarsely as follows: Multiply the Fourier transform of the given 
hyperfunction with an entire function which is then taken to be the reciprocal
of the Fourier transform of a local operator, which in turn is nothing but a 
multiplication operator with an entire function. This function can be chosen to
decay fast enough to ensure that its inverse Fourier transform is continuous on the 
real axis. To ensure the existence of these multipliers, we need two basic 
lemmata.
\begin{lemm}[{\cite[Proposition~8.1.6]{b:KAN88}}]
\labelL{K1}
Let $\phi(t)$ be a positive, monotonously increasing function on the
half-axis $t\geq1$ with $\phi(t)>1$ and
$\lim_{t\to\infty}\phi(t)=\infty$. 
Then the infinite product
\[
J(\zeta)=\prod_{k=1}^\infty
\left(1+\frac{\zeta^2}{(k\phi(k))^2}\right),
\]
with $\zeta^2\DEF\zeta_1^2+\ldots+\zeta_n^2$, 
is an infra-exponential, entire function which fulfills the lower
estimate
\[
\ABS{J(\zeta)}\geq
  C\exp\left(c\frac{\ABS{\zeta}}{\phi(\ABS{\zeta}+1)}\right)
\]
for $\ABS{\IM\zeta}\leq\max\{(1/\sqrt{3})\ABS{\RE\zeta},1\}$ and
constants $C$,~$c>0$.
\end{lemm}
\begin{lemm}[\cite{b:KAN88}, Lemma~8.1.7]
\labelL{K2}
Let $\{f_k(t)\}_{k\in\NN}$ a sequence of positive, continuous functions
on the half-axis $t\geq0$ that have infra-exponential growth.
Then exists a function $\phi(t)$ as in the assumption of \refL{K1}, and
constants $C_k$ such that
\[
f_k(t)\leq C_k\exp\left(\frac{t}{\phi(t+1)}\right)
\]
holds for $t\geq0$.
\end{lemm}
\begin{proof}[Proof of \refT{struct}]
We first consider the case $f\in\SQMI{\RR^n}$.
By decomposing $f$ linearly into components decreasing exponentially
outside chosen cones, we can assume that the
Fourier transform of $f$ can be represented as a single boundary
value
$
\widehat{f}(\xi)= \widehat{F}(\xi+i\Gamma 0),
$
and such that the holomorphic function  $\widehat{F}\in\HOIE{W}$ on
an infinitesimal wedge $W$ of type $\DD^n+i\Gamma$ is exponentially
decreasing outside another closed cone $\Delta^\circ$.
Let $\{K_j\}_{j\in\NN}$ be an exhausting sequence of compact sets
in $\QQ^n$ for $W$. Then, every one of the functions
$
h_j=\sup_{\eta\in K_j} \ABS{\smash{\widehat{F}(\xi+i\eta)}}
$
is infra-exponential and we can apply the two cited Lemmata to conclude
that there exists a positive, infra-exponential, entire function
 $J(\zeta)$ and constants $\{C_j\}_{j\in\NN}$ such that
$
\ABS{\smash{\widehat{F}_j(\xi+i\eta)}}\leq C_j\ABS{J(\xi+i\eta)}
$
holds for $\eta\in K_j$. By Theorem~\ref{theo:PW2},
$\widehat{F}(\xi)$ is a $C^\infty$-function on the real axis.
Consider the function
\[
\widehat{F}_{0}(\zeta)=
\frac{\widehat{F}(\zeta)}{J(\zeta)(1+\zeta^2)^n},
\]
which is for $\IM\zeta=\eta$ and $\eta\in W$ absolutely
integrable in $\xi=\RE\zeta$ and its restriction to
the real axis is $C^\infty$. Thus the inverse Fourier transform
of $\widehat{F}_{0}$ are well-defined in the sense of Fourier 
hyperfunctions as a boundary value $F_{0}=\FOUR^{-1}\widehat{F}_{0}$ 
on an infinitesimal wedge of type $\DD^n+i\Delta0$ which is furthermore
exponentially decreasing outside $-\Gamma_j^\circ$. 
Now, $F_{0}$ can be extended to a continuous function on the real axis,
which can be shown to be of asymptotic decay, since it is the
inverse Fourier transform of a smooth function. 
Set $f_0(x)=F_{0}(x)$, then Corollary~\ref{coro:fourlocop} shows
the claim. The case $f\in\SQPI{\RR^n}$ follows by similar reasoning, except
for the following modification: The Fourier transform $\FOUR f$ becomes
a distribution of finite order in this case and one can apply the well-known
structural theorem for these, see\eg~\cite{b:TRE67},
to represent it as a finite order differential 
operator applied to a continuous function. After dividing by
 $J(\zeta)(1+\zeta^2)^n$ and inverse Fourier transformation one can conclude
that $f_0$ becomes a continuous function with a certain polynomial\ie tempered,
growth.
\end{proof}
\section{Asymptotic Expansions}
\label{sect:asympt}
\subsection{One-Dimensional Asymptotic Expansions}
We are now ready to generalize the moment asymptotic expansions
of distributions~\cite{EK90,b:EK94} to the case of asymptotic hyperfunctions.
We start with the one-dimensional case.
\begin{defi}
For any asymptotic hyperfunction $f$ on \DD and every number $n=0,1,\ldots$, we define
its \textbf{$n$-th moment} by
\begin{displaymath}
  \MOM[f]{n}\DEF\int_\DD x^n\cdot f(x)\dd x
\end{displaymath}
(and write $\mu^n$ for short if there is no danger of confusion).
The space
\begin{displaymath}
  \SQ[{[N]}]{\DD}\DEF\Bigl\{ f\in\SQMI{\DD} \Bigm| \MOM[f]{n}=0,\
      n=0,\ldots,N-1\Bigr\}
\end{displaymath}
is called the \textbf{remainder space of order $N$}.
\end{defi}
Here, the integral over an asymptotic hyperfunction is defined as usual by the integral over a defining function. It is now an easy task to show the validity of asymptotic expansions for hyperfunctions on \DD.
\begin{theo}
\label{theo:asymptexp}
The \textbf{moment asymptotic expansion}
\begin{gather*}
  f(x)=\SUM[f]{N}(x) \mod \SQ[{[N+1]}]{\DD},\\
  \intertext{with the \textbf{asymptotic sum of order $N$}, given by}
  \SUM[f]{N}(x)\DEF\sum_{n=0}^N
  \frac{(-1)^n}{n!}\MOM[f]{n}\delta^{(n)}(x),
\end{gather*}
holds for any $f\in\SQMI{\DD}$.
\end{theo}
\begin{proof}
We have to show that the remainder of the $N$-th order expansion is
in the remainder space of the same order:
\begin{displaymath}
  \REM[f]{N}(x)\DEF f(x)-\SUM[f]{N}(x)\in\SQ[{[N+1]}]{\DD},
\end{displaymath}
for every $N\in\NN$. The partial sum \SUM[f]{N} is trivially a hyperfunction
with support in the origin. Since the exponentially decaying Cauchy--Hilbert kernel $h_z$ of the last section induces a continuous embedding of $\SB[\ast]{\RR}$ into \SQMI{\DD}, we can regard \SUM[f]{N} as an asymptotic hyperfunction. Thus the difference $\REM[f]{N}=f-\SUM[f]{N}$ is again asymptotic and $\REM[f]{N}\in\SQ[{[N+1]}]{\DD}$ follows immediately from the definition of the remainder.
\end{proof}
In applications of asymptotic expansions, as well as for the expansions of distributions as in~\cite{EK90,b:EK94,EGV98}, it is useful to state the expansion in the so-called \emph{parametric form}\ie in application to a test function with scaled argument. We start with a technical lemma:
\begin{lemm}
For any continuous seminorm $\NORM{\cdot}$ on $\SPPI{\RR}$
and every function $\phi\in\SPPI{\RR}$ with $\phi^{(n)}(0)=0$ for $n=0,\ldots,N$ holds
\[
\NORM{\phi(z/\lambda)}=O(\ABS{\lambda}^{-N-1}),
\]
for $\ABS{\lambda}\to\infty$.
\end{lemm}
\begin{proof}
Similar to the proof of Theorem~\ref{theo:PW1}, we represent
\SPPI{\DD} as an inductive limit $\SPPI{\DD}\cong\varinjlim
\smash[t]{\HOBO[m]{\ICLOS{U_m}}}$, $m\in\NN$, with neighborhoods
$U_m=\{z\in\CC\mid\ABS{\IM z}<1/m\}$. A locally convex topology
on $\smash[t]{\HOBO[m]{\ICLOS{U_m}}}$ can be generated by the system of seminorms
\[
\NODX{\phi(z)}{m}{R}\DEF \sup\left\{\ABS{\smash{\phi(z)(1+\ABS{\RE
z}^{-m})}}\Bigm| z\in \overline{U_m},\ \ABS{\RE z}\leq R\right\},
\]
for $R>0$. This topology is apparently weaker than the original
topology of the space \HOBO[m]{\ICLOS{U_m}}, but has the following property,
which is sufficient for our purpose: If $\phi_n$ is a sequence
which is bounded in \HOBO[m]{\ICLOS{U_m}} and converges to zero
in the seminorms \NODX{.}{m}{R} for every $R>0$, then it converges 
to zero in \HOBO[m']{\ICLOS{U_{m'}}} for every $m'>m$,
cf.\ Definition~\ref{defi:spaces}. That means,
the systems of seminorms \NODX{.}{m}{R} induces on the bounded 
subsets of \SPPI{\DD} the original topology of that space.
Now, for any $\phi\in\SPPI{\DD}$ the set $\{\phi(z/\lambda)\}_{\lambda>C>0}$
is bounded in \SPPI{\DD} and therefore it is also bounded in some
\HOBO[m]{\ICLOS{U_m}}. We can thus argue as follows:
Let $\phi\in\HOBO[m]{\ICLOS{U_m}}$ such that $\phi^{(n)}(0)=0$ for 
$n=0,\ldots,N$. Then, there is a constant $K>0$ for which
\[
\ABS{\phi(z)}\leq K\ABS{\RE z}^{N+1},\ z\in U_m,\ \ABS{\RE z}\leq 1.
\]
If $\lambda>R$ then $\NODX{\phi(z/\lambda)}{m}{R}\leq
K/\lambda^{N+1}$, and thus
\[
\NODX{\phi(z/\lambda)}{m}{R}=O(\ABS{\smash{\lambda^{-N-1}}})
\]
for  all $R$. Since the topology of $\SPPI{\RR}$ is generated by all seminorms 
for which the restriction to every subspace \HOBO[m]{\ICLOS{U_m}} is again a 
continuous seminorm, the assertion follows.
\end{proof}
A lemma of this type is fundamental for every asymptotic expansion in parametric form, as the following general arguments will show.
\begin{theo}
\label{theo:Parameterform}
For every $f\in\SQMI{\DD}$, the moment asymptotic expansion in
\textbf{parametric form}
\[
f(\lambda x)\sim\sum_{n=0}^\infty
\frac{(-1)^n\MOM[f]{n}\delta^{(n)}(x)}{n!\lambda^{n+1}},
\]
holds for $\lambda\in\RR$, $\ABS{\lambda}\to\infty$.
This formula holds in the dual sense\ie for all
$\phi\in\SPPI{\DD}$ holds
\[
\SPROD{f(\lambda x)}{\phi(x)}=
\sum_{n=0}^N \frac{\MOM[f]{n}\phi^{(n)}(0)}{n!\lambda^{n+1}}+
O\left(\ABS{\lambda}^{-N-2}\right),
\]
for $\lambda\in\RR$, $\ABS{\lambda}\to\infty$.
\end{theo}
\begin{proof}
If $f$ and $\phi$ are as above then
\[
\SPROD{f(\lambda x)}{\phi(x)}=\frac{1}{\lambda}
\SPROD{f(x)}{\phi(\lambda^{-1}x)}.
\]
We write the Taylor expansion of $\phi(x/\lambda)$ in $x$ around $0$ with remainder as
\[
\phi(x/\lambda)=
\sum_{n=0}^N \frac{\phi^{(n)}(0)}{n!\lambda^n}\cdot x^n+
T_N(\lambda;x).
\]
Since $\REM[f]{N}\in\SQ[{[N+1]}]{\DD}$, it follows
\begin{align*}
\SPROD{f(\lambda x)}{\phi(x)}&=
\SPROD{\SUM[f]{N}(\lambda x)+\REM[f]{N}(\lambda x)}{\phi(x)}\\
&=
\sum_{n=0}^N \frac{\MOM[f]{n}\phi^{(n)}(0)}{n!\lambda^{n+1}}+
\lambda^{-1}\SPROD{\REM[f]{N}(x)}{T_N(\lambda;x)}.
\end{align*}
Now, $T_N(\lambda;x)$ is in $\SPPI{\DD}$ and $D^n_x
T_N(\lambda;0)=0$, $n=0,\ldots, N$. Consequently, the assertion
follows from the above lemma, the definition of $T_N$ and the
fact that \REM[f]{N} is a continuous linear functional on
\SPPI{\DD}.
\end{proof}
We state clearly at this point that, despite the synonymy, our
asymptotic expansions of hyperfunctions have almost nothing to
do with that based on second microlocalization, see~\cite{KK80}.
Rather it is a proper generalization of the distributional expansions by
R.~Estrada and coworkers~\cite{EK90,b:EK94}.

It is simple but useful to restate the expansion for the Fourier transforms
of asymptotic hyperfunctions.
\begin{prop}
\label{prop:Fourasympt}
Let $f\in\SQMI{\DD}$. Then, its Fourier transform $\widehat{f}$
fulfills the asymptotic expansion
\begin{displaymath}
  \widehat{f}(\xi)=\sum_{n=0}^{N}
  \frac{(-i\xi)^n}{n!}\MOM[f]{n}\mod
  \FQIDX{[N+1]}{\RR}
\end{displaymath}
around $0$\ie for $\ABS{x}\to0$. Here, we denote
the image of \SQ[{[N+1]}]{\DD} under the Fourier transformation
by
\begin{displaymath}
  \FQIDX{[N+1]}{\RR}\DEF\Bigl\{
    g\in\CIIE{\RR}\Bigm| g^{(n)}(0)=0,\ n=0,\ldots,N
\Bigr\}.
\end{displaymath}
\end{prop}
\begin{proof}
By Theorem~\ref{theo:PW2}, we can calculate the derivatives
of $\widehat{f}$ at $0$ explicitly:
\[
  i^k\widehat{f}^{(k)}(0)=
    \left(i^k\frac{\dd^k}{\dd\xi^k}\int_\DD \ee^{-ix\xi}f(x)\dd x
    \right)\Bigg|_{\xi=0}
  = i^k\int_\DD (-ix)^k f(x) \dd x = \MOM[f]{k}.
\]
That is, the stated asymptotic expansion is identical to the
Taylor expansion of the smooth function $\widehat{f}$ around $0$ and we have
$\widehat{R}_{N}^{(n)}(0)=0$ for $n\leq N$ by definition of the remainder
space \SQ[{[N+1]}]{\DD}.
\end{proof}
This allows us to derive the following basic result, cf.~\cite{EK90,b:EK94}:
\begin{prop}
\label{prop:complete}
The  moment asymptotic expansion of asymptotic hyperfunctions
is a \textbf{complete asymptotic scheme}\ie for any given sequence
$\{\MOM{n}\}_{n\in\NN_0}\subset\CC$ exists an
$f\in\SQMI{\DD}$ with $\MOM[f]{n}=\MOM{n}$.
\end{prop}
In order to prove this assertion, we use the classical Theorem of
Ritt in a formulation which can be derived from~\cite[Section
1]{DE94}.
\begin{theo}
\label{theo:temperedmoment}
For $\{\MOM{n}\}_{n\in\NN_0}\subset\CC$ and a given sector
$S=\{\zeta\in\CC\mid \alpha<\arg\zeta<\beta,\ 0<\ABS{\zeta}\}$
with vertex $0$ in the complex plane, there
exists a function $\psi(\zeta)$, which is holomorphic
in $S$, bounded in $\overline{S}$, satisfies
$\psi(\zeta)=O(\ABS{\zeta}^{-\infty})$
for $\zeta\to\infty$ in $S$, and such that furthermore holds
\[
\psi(\zeta)\sim\sum_{n=0}^{\infty}\frac{(-i\zeta)^n}{n!}\MOM{n},
\]
for $\zeta\to 0$ in $S$.
\end{theo}
\begin{proof}[Proof of Proposition~\ref{prop:complete}]
By Proposition~\ref{prop:Fourasympt} and Theorem~\ref{theo:PW2}, the
assertion is equivalent to: There exists a function
$\widehat{f}\in\CIIE{\RR}$ such that
$\widehat{f}^{(n)}(0)=(-i)^n\MOM{n}$, since then
$f=\FOUR^{-1}\widehat{f}\in\SQMI{\DD}$ fulfills the original condition.
Choose the sector $S$ in Theorem~\ref{theo:temperedmoment} large enough
to contain $\RR\setminus\{0\}$ and a function $\psi$ as there.
This function is analytic and therefore $C^\infty$ in $\RR\setminus\{0\}$
and has a $C^\infty$-continuation to the point $0$.
Furthermore, it is asymptotic along the real axis, and thus
we can choose $\widehat{f}=\psi\in\CIIE{\RR}$.
\end{proof}
We remark that $f$ can even be chosen to lie in $\SSS{\RR}$.

\subsection{Radon Transformation of Asymptotic Hyperfunctions}
Our final task is to carry over the one-dimensional asymptotic expansions
to higher dimensions. To this end, we use the geometrical
\textbf{Radon transformation} of hyperfunctions, introduced by Kaneko and
Takeguchi, see~\cite{KT95}. This transformation decomposes a given
asymptotic hyperfunction into a one-dimensional part depending only on the
radial coordinate and a second part depending on the remaining angle
coordinates. For the readers convenience, we briefly review the general concept.
Note that the following results, up to Definition and Theorem~\ref{theo:Radon} are proved
in~\cite{KT95}.

For any $f\in\SSS{\RR^n}$, its \textbf{Radon transform}
is defined by integrals over the affine planes
$P(\omega,t)\DEF \bigl\{x\in\RR^n\bigm| \omega x=t\bigr\}$
for $\omega\in S^{n-1}$ and $t\in\RR$ by
\begin{displaymath}
  \RAD f(\omega,t)\DEF\int_{\RR^n}\delta(t-\omega x)f(x) \dd x.
\end{displaymath}
$\RAD f$ is a $C^\infty$-function on $S^{n-1}$ with values in
\SSS{\RR}\ie an element of the topological vector space
$C^\infty(S^{n-1};\SSS{\RR})$ in the notation of~\cite[Chapters
40--44]{b:TRE67}. Since \SSS{\RR} is nuclear, this space is
isomorphic to a completed tensor product and thus $\RAD f\in
C^\infty(S^{n-1})\CTENS\SSS{\RR}$ holds. It is geometrically
clear that $\RAD f$ is an even function in $(\omega,t)$\ie $\RAD
f(-\omega,-t)=\RAD f(\omega,t)$.

For a hyperfunction $f$, neither its restriction to $P(\omega,t)$
nor its integral over this hyperplane is well defined in general,
see~\cite[Chapter 3, §4]{b:KAN88}. To remedy that problem, Kaneko
and Takiguchi introduce a special class of \emph{Radon
hyperfunctions}, for which they define the Radon transformation
via duality. These Radon hyperfunctions are not a proper subclass
of Fourier hyperfunctions, but for appropriate subclasses of
Radon hyperfunctions which are also Fourier hyperfunctions, the
Radon transformation can be defined explicitly using the
boundary value representation of hyperfunctions. This is
especially the case for asymptotic hyperfunctions, as we will now
see.

Using the decomposition
\begin{displaymath}
  \delta(t-\omega x)=\MTPI\left(
    \frac{1}{(t+i0)-\omega x} - \frac{1}{(t-i0)-\omega x}\right)
\end{displaymath}
of the $\delta$-function, we are led to the formal definition
of the Radon transform for an $f\in\SQMI{\DD^n}$:
\begin{align*}
  \RAD f(\omega,t)&\DEF G(\omega,t+i0)-G(\omega,t-i0),\\
\intertext{with}
  G(\omega,\tau)&\DEF\MTPI\int_{\DD^n}\frac{f(x)}{\tau-\omega x}\dd x,
\end{align*}
for $\tau\in\CC\setminus\RR$. This result is to be understood as a Fourier hyperfunction
valued ordinary function in $\omega$.
In a special boundary value representation $f(x)=\sum_{j=1}^N
F_j(x+i\Gamma_j0)$, with $F_j\in\HOMI{\DD^n+i\Gamma_j0}$, this becomes
\begin{gather*}
  \RAD f(\omega,t)\DEF \sum_{j=1}^N G_j(\omega,t+i0)-G_j(\omega,t-i0),\\
\intertext{with defining functions $G_j$ given by}
  G_j(\omega,\tau_j)
   \DEF\MTPI\int\limits_{\IM z_j=y_j}\frac{F_j(z_j)}{\tau_j-\omega z_j}\dd z_j
   =\MTPI\int\limits_{\RR^n} \frac{F_j(x+iy_j)}{t+is_j-\omega (x+iy_j)}\dd x.
\end{gather*}
Here, we set $\tau_j=t+is_j$ and $z_j=x+iy_j$ for $y_j\in\Gamma_j$.
The integrals converge absolutely and independently of the damping factor
$(\tau_j-\omega z_j)^{-1}$, and yield a defining function for $\RAD f$.
Furthermore, the definition is independent of the choice of integration planes
${\IM z_j=y_j}$ by Cauchy's Theorem.
In the sense that for $\omega\to\omega_0\in S^{n-1}$ one has
$G_j(\omega,\tau_j)\to G_j(\omega_0,\tau_j)$ in the topology
of \HOIE{\QQ\setminus\DD}, we can regard $\RAD f(\omega,t)$ as a
continuous function on $S^{n-1}$. That is, we have
$\RAD f\in C(S^{n-1};\SQ[]{\DD_t})\cong C(S^{n-1})\CTENS\SQ[]{\DD_t}$,
where $\DD_t$ denotes \DD with coordinate $t$. We will see below that
the Radon transformation is well-defined\ie does not depend on the choice
of boundary value representation of $f$, by exhibiting its connection
with the Fourier transformation.

To see that $\RAD f$ is an asymptotic hyperfunction in $t$ for all
$\omega$, we first recover the original representation of the Radon
transform by integrals over hyperplanes from the above definition of the $G_j$:
\begin{align*}
  G_j(\omega,\tau_j)&=
  \int\limits_{\RR^n}\delta(t+is_j-\omega(x+iy_j))F_j(x+iy_j)\dd x\\&=
  \int\limits_{P(\omega,\tau_j)} F_j(x+iy_j)\dd S.
\end{align*}
Here, $P(\omega,\tau_j)$ is the $n-1$-dimensional real affine plane
\[
P(\omega,\tau_j)\DEF \bigl\{z_j=x+iy_j\bigm|
t+is_j=\omega(x+iy_j)\bigr\}\subset\CC^n.
\]
Now, the estimate $\ABS{x}=O(\ABS{t})$ holds locally uniformly in $s_j$
and $y_j$ on $P(\omega,\tau_j)$. From this it follows, that
$G_j$ is an asymptotic function in $t$ locally uniformly in $y_j$.
Thus indeed, $\RAD f(\omega,\cdot)\in\SQMI{\DD_t}$
for all $\omega\in S^{n-1}$.

There is a close connection between the Radon and the Fourier
transformations, which is given by the formula
\begin{displaymath}
  \RAD f(\omega,t)=\frac{1}{2\pi}\int_{-\infty}^{\infty}
  \widehat{f}(\rho\omega)\ee^{it\rho}\dd\rho.
\end{displaymath}
For $f(x)\in\SQMI{\DD^n}$ we have $\widehat{f}(\xi)\in\CIIE{\RR^n}$
by Theorem~\ref{theo:PW2}, and its partial back transformation
along the fiber $\xi=\rho\omega$ can be interpreted as a Fourier
hyperfunction as follows: Split the integration path at $0$ and
compute the two components
\begin{displaymath}
  G_±(\omega,\tau)\DEF±\frac{1}{2\pi}\int_0^{±\infty}
  \widehat{f}(\rho\omega)\ee^{i\tau\rho}\dd\rho
\end{displaymath}
of the defining function of $\RAD f$ for $±\IM\tau>0$ separately.
We show that this yields the same result as the original formula for
$G(\omega,\tau)$:
\begin{align*}
G_±(\omega,\tau)&=
\frac{1}{2\pi}\int_0^{±\infty}\ee^{i\tau\rho}
  \left( \int_{\DD^n} f(x) \ee^{i\rho\omega x} \dd x
  \right) \dd\rho, \\
\intertext{for $±\IM\tau>0$. The inner integral is an absolutely
convergent integral over defining functions and can thus be exchanged with
the outer one, which leads to}
G_±(\omega,\tau)&=\int_{\DD^n}f(x)\left(
\frac{1}{2\pi}\int_0^{±\infty}\ee^{i\rho(\tau-\omega
x)} \dd\rho\right)\dd x\\
&=\int_{\DD^n} \MTPI\frac{f(x)}{\tau-\omega x}\dd x=G(\omega,\tau),
\end{align*}
as we wanted to show. As an aside, this also shows the invariance of
the definition of $\RAD f$ under a change of the boundary value
representation of $f$ by the namely property of the Fourier transformation,
see\eg~\cite[Lemma~8.3.3 and Theorem~8.3.4]{b:KAN88}.
We compile the information we won so far in the following
\begin{defitheo}\label{theo:Radon}
We define the \textbf{Radon transform} of $f(x)\in\SQMI{\DD^n}$, denoted by
$\RAD f(\omega,\cdot)$, by the functions
$G_±(\omega,\tau_±)\in C^\infty(S^{n-1};\HOMI{\QQ_±})$,
called \textbf{canonical defining functions},
and where $\QQ_±\DEF\ICLOS{\{z\mid\IM{z}\gtrless0\}}\setminus\DD$,
via the boundary values
\begin{gather*}
\RAD f(\omega,t)=[G(\omega,\tau)]_{\tau=t}=
[G_+(\omega,t+is),G_-(\omega,t-is)]_{s=0}.
\end{gather*}
We have
$\RAD f\in C^k(S^{n-1})\CTENS\SQMI{\DD_t}$,
$\RAD f$ is an even function in $(\omega,t)$,
and the canonical defining function $G(\omega,\tau)$
fulfills the estimates
\begin{displaymath}
  \ABS{D_\omega^\alpha G(\omega,\tau)}\leq
  C\frac{(\ABS{\alpha}!)^2}{\nu^{\ABS{\alpha}}},
\end{displaymath}
with constants $\nu$, $C>0$.
\end{defitheo}
\begin{proof}
We only have to show the bounds on the derivatives of $G$. The
proof follows~\cite[Propositions 2.3 and 2.8]{KT95}. As in
Theorem~\ref{theo:PW2}, we can assume $f$ to be represented by a
single boundary value $f=F(x+i\Gamma 0)$,
$F\in\HOMI{\DD^n+i\Gamma0}$, such that $F$ decreases
exponentially outside a certain cone $\Delta^\circ$. A coarse
estimate for the Fourier transform of $F$ is
\begin{align*}
  \ABS{\smash{D^\alpha_\xi \widehat{F}}(\zeta)}&\leq
  \ABS{\int_{\RR^n} \ee^{-iz\zeta}(-iz)^\alpha F(z) \dd x} \leq
  C\int_{\RR^n}\ee^{x\eta+y\xi}\ABS{x^\alpha}\ABS{F(x+iy)}\dd x,\\
\intertext{with $\zeta=\xi+i\eta$ and  $z=x+iy$. Now, for every compact set
$L\subset-\Delta$ and $\eta\in L$, there is a constant $\delta_L>0$ such that
we can further estimate}
  \ABS{\smash{D^\alpha_\xi \widehat{F}}(\zeta)}&
  \leq C'\ee^{y\xi}\int_{\RR^n}\ee^{-\delta_L\ABS{x}}\ABS{x^\alpha}
  \ABS{F(x+iy)}\dd x
  \leq C''\cdot\frac{\ABS{\alpha}!}{\delta_L^{\ABS{\alpha}}}\cdot\ee^{y\xi},
\end{align*}
cf.\ the proof of~\cite[Theorem 8.3.2]{b:KAN88}. Since the
derivatives $D^\alpha_\omega\widehat{f}(\rho\omega)$ exist for
arbitrary $\alpha$, we can calculate that of\eg $G_+$:
\begin{align*}
  \ABS{D_\omega^\alpha G_+(\omega,t+is)}&=
  \ABS{\frac{1}{2\pi}D_\omega^\alpha \int_0^\infty \ee^{i(\rho+i\sigma)(t+is)}
    \widehat{F}((\rho+i\sigma)\omega) \dd\rho},\\
\intertext{where $\omega\in-\Delta$ if we choose $\sigma>0$. This can be
estimated by}
&\leq
  \ABS{\frac{1}{2\pi}\int_0^\infty \ee^{-\rho s-\sigma t}
  \rho^{\ABS{\alpha}}
  (D_\xi^\alpha\widehat{F})((\rho+i\sigma)\omega)\dd\rho}.\\
\intertext{Using the first estimate above, we further get for every $\kappa>0$}
  &\leq C\int_0^\infty \rho^{\ABS{\alpha}}
  \frac{\ABS{\alpha}!}{\delta_\sigma^{\ABS{\alpha}}} \ee^{\kappa\rho}
  \ee^{-\rho s-\sigma t}\dd\rho
  \leq C'
  \frac{(\ABS{\alpha}!)^2}{(\delta_\sigma(s-\kappa))^{\ABS{\alpha}}}
    \ee^{-\sigma t}.
\end{align*}
Since this final result is finite for every $\sigma$, we can replace
the denominator by $\nu^{\ABS{\alpha}}$ for a suitable constant $\nu>0$.
This proves the assertion.
\end{proof}
The regularity of $\RAD f$ expressed by the bounds on its derivatives in $\omega$, is
the well known one of \textbf{ultradifferentiable functions}, cf.~\cite{KOM89}:
\begin{defi}
 A function $g\in C^\infty(\RR^n)$ is said to be in the \textbf{Gevrey-Class of
order $s$} for $s>1$, if and only if for every compact set $K\subset\RR^n$ there are
constants $h$ and $C>0$ such that
\begin{displaymath}
  \NORM{D^\alpha g}_{C^0(K)}\leq Ch^{\ABS{\alpha}}(\ABS{\alpha}!)^s
\end{displaymath}
holds for $\ABS{\alpha}=0,1,2,\ldots$.
\end{defi}
For $f\in\SQMI{\DD^n}$, we immediately find by Theorem~\ref{theo:Radon}
that Gevrey-Bounds hold for all derivatives of $\RAD
f(\omega,\cdot)$ with respect to $\omega\in S^{n-1}$, for $h=\nu^{-1}$ and
with $s=2$, that is:
\begin{coro}
\label{coro:QMIRadon}
For $f\in\SQMI{\DD^n}$, $\RAD f(\omega,t)$ is a function in the Gevrey-Class
of order $2$ on $S^{n-1}$ with values in $\SQMI{\DD_t}$.
\end{coro}
This result is analogous to~\cite[Proposition 2.8]{KT95}, where
it is shown that this type of regularity holds for the Radon
transforms of exponentially decaying hyperfunctions, and is the
best possible in this case.
\subsection{Radon Asymptotic Expansions}
The moments of a Radon transform $\RAD f$, which are the last ingredients
for asymptotic expansions in higher dimensions are known as the \textbf{Helgason
moments}, see~\cite{b:HEL84}, and are defined by the formula
\begin{gather*}
  \HELMOM{f}{k}\DEF \MOM[\RAD f(\omega,\cdot)]{k}
  =\int_{\DD_t} t^k\cdot
  \RAD f(\omega,t)\dd t,\quad\text{for }k=0,1,2,\ldots.
\end{gather*}
The \HELMOM{f}{k} satisfy the \textbf{Helgason moment condition},
cf.~\cite[p.\ 100]{b:HEL84}:
\begin{prop}
\label{prop:Helgason}
For $f\in\SQMI{\DD^n}$ and $k=0,1,2,\ldots$, the Helgason moment \HELMOM{f}{k} is
a homogeneous polynomial of total degree $k$ in $\omega$.
\end{prop}
\begin{proof}
We first calculate the Helgason moments as in the proof of
Proposition~\ref{prop:Fourasympt}, to obtain
\begin{align*}
 \HELMOM{f}{k}&=
  \left(\int_{\DD_t} \RAD f(\omega,t) t^k\ee^{-it\rho}\dd t
  \right)\biggr|_{\rho=0} =
  \left(i^k\frac{\dd^k}{\dd\rho^k}\int_{\DD_t} \RAD f(\omega,t)\ee^{-it\rho}\dd t
  \right)\Biggr|_{\rho=0}\\
  &= \left(i^k\frac{\dd^k}{\dd\rho^k}\widehat{f}(\rho\omega)
  \right)\Biggr|_{\rho=0}.
\end{align*}
Since $\widehat{f}$ is $C^\infty$,
we can write down its Taylor series up to the order $k$ with remainder
around the point $0$ as
\begin{displaymath}
  \widehat{f}(\xi)=\sum_{\ABS{\alpha}\leq k} a_\alpha
  \xi^\alpha+R_{k+1}(\xi),
\end{displaymath}
with $a_\alpha=D_\xi^\alpha \widehat{f}(0)/\alpha!$ and
\begin{displaymath}
  R_{k+1}(\xi)=\sum_{\ABS{\beta}= k+1} \frac{D_\xi^\beta
    \widehat{f}(\Theta\xi)}{\beta!} \xi^\beta,
\end{displaymath}
for a suitable $\Theta\in[0,1]$. Inserting this in the formula above yields,
since we have $D^k_\rho R_{k+1}(\rho\omega)=O(\ABS{\rho})$ for $\ABS{\rho}\to 0$,
the exact form of \HELMOM{f}{k}, namely:
\[
\HELMOM{f}{k}=
i^kk!\sum_{\ABS{\alpha}=k}a_\alpha\omega^\alpha,
\]
proving the assertion.
\end{proof}
Now, we have collected all necessary facts about the Radon transformation of
asymptotic hyperfunctions, their Helgason moments, and the relation to
their Fourier transforms to carry over
the asymptotic expansions of Theorem~\ref{theo:asymptexp} and
Proposition~\ref{prop:Fourasympt} to the case of dimension $n>1$.
\begin{theo}
  For $f\in\SQMI{\DD^n}$, the \textbf{Radon asymptotic expansion}
\[
  \RAD f(\omega,t) = \sum_{k=0}^N \frac{(-1)^k}{k!} \HELMOM{f}{k} \cdot
    \delta^{(k)}(t) \mod \RAD \SQ[{[N+1]}]{S^{n-1}×\DD_t}
\]
holds good.
Here, we defined the \textbf{Radon remainder space of order} $N+1$ by
\begin{displaymath}
    \RAD \SQ[{[N+1]}]{S^{n-1}×\DD_t}\DEF \left\{
    \RAD g\in\RAD\SQMI{\DD^n}\Bigm| \HELMOM[]{g}{k}\equiv0,\
    k=0,\ldots,N\right\}.
\end{displaymath}
Furthermore, the following relation to the Taylor expansion of
the Fourier transform
$\widehat{f}$ around the origin holds:
\[
\RAD f(\omega,t) = \sum_{k=0}^N \sum_{\ABS{\alpha}=k} (-i)^k a_\alpha
     \omega^\alpha \cdot \delta^{(k)}(t)
     \mod \RAD \SQ[{[N+1]}]{S^{n-1}×\DD_t},
\]
where $a_\alpha=D_\xi^\alpha \widehat{f}(0)/\alpha!$.
\end{theo}
\begin{proof}
By~\cite[Theorem 2.4]{KT95}, the $N$-th partial sum
$\smash{\SUM[\RAD f]{N}}$ of the Radon asymp­to­tic expansion
has a unique preimage under the Radon transformation. By the
same theorem, this preimage $\RAD^{-1}\smash{\SUM[\RAD f]{N}}$ is
a hyperfunction with compact support and therefore especially in
\SQMI{\DD^n}. Thus, the remainder $f-\RAD^{-1}\SUM[\RAD f]{N}$ is
asymptotic and its image under $\RAD$ is exactly the remainder of
the Radon asymptotic expansion. The Helgason moments of the
remainder vanish up to order $N$ by definition.
\end{proof}
Note aside how every partial sum $\smash{\SUM[\RAD f]{N}}$ becomes an even function
in $(\omega,t)$: When changing the variable to
$(-\omega,-t)$, the sign of $\delta^{(k)}$ cancels that of $\HELMOM[]{f}{k}$.
We refrain from rewriting the Radon asymptotics in parametric form as
in Theo­rem~\ref{theo:Parameterform}, and instead give a very basic example.
\begin{nxmpl}
Let $f(x)=J(D)\delta(x-a)$, with $a\in\RR^n$, and let $J(D)=\sum_{\ABS{\alpha}\geq0}b_\alpha D_x^\alpha$
be a local operator. The Radon transform of $f$ is
\[
\RAD f(\omega,t)=J(\omega D_t)\delta(t-a\omega),
\]
see~\cite[Example 3.3]{KT95}, by which the Helgason moments can
easily be calculated. This yields the Radon asymptotic expansion
\[
\RAD f(\omega,t)\sim \sum_{k=0}^\infty \sum_{\ABS{\alpha}\leq k}
\frac{(-1)^{k-\ABS{\alpha}}}{\ABS{\alpha}!}\cdot
(b_\alpha\omega^\alpha)(a\omega)^{k-\ABS{\alpha}}\cdot\delta^{(k)}(t).
\]
This representation makes the Helgason moment condition manifest.
\end{nxmpl}
We conclude our analysis by giving a condition on the support of
a hyperfunction in terms of the moments of its Radon asymptotic expansion.
To do so, we apply a recent result by Kim,~\emph{et al.}:
\begin{theo}[{\cite[Theorem 3.1]{CKY98}}]
\label{theo:CKY}
A sequence $\{\MOM{k}\}_{k\in\NN_0}$ is the  moment sequence
of a hyperfunction $f\in\SB[\ast]{\RR}$ with support in the interval
$[-R,R]$, if and
only if for every $S>R$ and $\EPS>0$ there exists a constant $C_\EPS>0$
such that the estimate
\begin{gather*}
\ABS{\sum_{k=0}^\infty\frac{\MOM{k}}{k!}\frac{1}{2S}
\left(-\frac{\pi iq}{S}\right)^k}\leq
C_{\EPS} \ee^{\EPS\ABS{q}}
\tag{\ddag}\label{eq:momest}
\end{gather*}
holds for all $q\in\ZZ$.
\end{theo}
This result was originally formulated for cubic domains in
$\RR^n$, but for our purpose a one dimensional version suffices,
since we want to use the Radon transformation again to generalize
to higher dimensions. To that end, we need an adaption of another
result of~\cite{KT95}, which connects the support of a Radon
transform with that of its Radon--preimage. This so called
\textbf{support theorem} goes back to Helgason, see~\cite[pp.\
105]{b:HEL84}, where it is proved for functions in \SD{\RR^n}.
\begin{theo}[{\cite[Theorem 4.1]{KT95}}]
\label{theo:support}
If $f\in\SQMI{\DD^n}$ is such that $\RAD f(\omega,t)$ vanishes for
$\ABS{t}\geq R$, then $f(x)$ vanishes for $\ABS{x}\geq R$.
\end{theo}
We finally give the intended result which allows one to restrict the support
of hyperfunctions to spherical domains, provided that a growth condition
on the Helgason moments of the Radon asymptotic expansion holds:
\begin{theo}
\label{theo:Radonsuppcond}
A sequence $\{\alpha^k(\omega)\}_{k\in\NN_0}$ of polynomials on
$S^{n-1}$ which fulfill Helgason's moment condition is the
sequence of Helgason moments of the Radon transform $\RAD f$
of a hyperfunction $f\in\SB[\ast]{\RR^n}$ with support in the ball
$B_R(0)=\{x\in\RR^n\mid\ABS{x}\leq R\}$,
if and only if the $\alpha^k$ satisfy the estimate~\eqref{eq:momest}.
\end{theo}
\begin{proof}
If $f\in\SQMI{\DD^n}$  has support in $B_R(0)$ then $\RAD
f(\omega,t)$ vanishes for $\ABS{t}>R$ and all $\omega$.
The estimate is then immediate from Theorem~\ref{theo:CKY}.
Let conversely $\alpha^k(\omega)$ be a sequence of polynomials satisfying
Helgason's condition and the estimate~\eqref{eq:momest}.
From Proposition~\ref{prop:complete} and the formula for the Helgason moments
found in the proof of Proposition~\ref{prop:Helgason} it follows that
there is an asymptotic hyperfunction $g\in\SQMI{\DD^n}$
such that $\HELMOM{g}{k}=\alpha^k$. By the assumption and Theorem~\ref{theo:CKY},
$g$ can be chosen such that every component $\RAD g(\omega,\cdot)$ has compact
support in $[-R,R]$. Then $g$ itself vanishes for
$\ABS{x}>R$ by Theorem~\ref{theo:support}.
\end{proof}
We finally want to give a very simple example that shows how the
Radon asymptotic expansion can be utilized to solve differential equations
by an asymptotic series ansatz. Naturally, the Radon transform becomes most
effective when the equation in question exhibits a spherical
symmetry.
\begin{nxmpl}
Consider the ordinary differential equation on $\DD^n$
\[
\left(r^2\frac{\dd}{\dd r}-1\right)f(x)=0.\quad\text{where }
r=\left(\sum_{i=1}^n x_i^2\right)^{1/2}.
\]
After Radon transformation this is easily seen to correspond to the equation
\[
\left(t^2\frac{\dd}{\dd t}-1\right)\RAD f(\omega,t)=0
\]
on $S^{n-1}×\DD_t$ only depending on $t$\ie an one-dimensional
equation. It is well-known, see~\cite[Example~3.9.7]{b:MOR93}, that
this irregular-singular equation has two pure hyperfunction solutions
which we try to recover by an asymptotic series ansatz.
We first calculate
\begin{align*}
t^2 \frac{\dd}{\dd t}\delta^{(n)}(t)&=
\left[\MTPI\frac{\tau^2(-1)^{n+1}(n+1)!}{\tau^{n+1}}\right]
= (n+1)n\left[\MTPI\frac{(-1)^{n-1}(n-1)!}{\tau^{n-1}}\right]
\\& =
\begin{cases}
0 & \text{ for } n\leq 1;\\ (n+1)n\cdot\delta^{(n-1)}(t) &\text{
otherwise,}
\end{cases}
\end{align*}
where we used the usual notation for hyperfunctions in terms of defining functions.
With that the equation becomes
\begin{align*}
\left(t^2\frac{\dd}{\dd t}-1\right)\sum_{n=0}^\infty d_n\delta^{(n)}(t)&=
\sum_{n=1}^\infty d_n(n+1)n\delta^{(n-1)}(t)-
\sum_{n=0}^\infty d_n\delta^{(n)}(t)\\
&=
\sum_{n=1}^\infty \left(d_n(n+1)n-d_{n-1}\right)\delta^{(n-1)}(t)=0,
\end{align*}
yielding the recursive prescription $d_n=d_{n-1}/(n(n+1))$ for the
coefficients. We solve it with initial condition $d_0=1$ through setting
$d_n=((n+1)!n!)^{-1}$, for all $n=0,1,\ldots$.
These coefficients decay fast enough to turn 
$J=\sum_{n=0}^\infty d_n D^n_t$ into a local operator and thus
the asymptotic series in this case actually converges to
a hyperfunction with support in the single point $\{0\}$.
This yields the first hyperfunction solution:
\begin{align*}
f_1(\omega,t)&=\sum_{n=0}^\infty d_n\delta^{(n)}(t)=
\frac{1}{2\pi
i}\left[\sum_{n=0}^\infty\frac{(-1)^{n+1}}{(n+1)!\tau^{n+1}}\right]_{\tau=t}\\&=
\frac{1}{2\pi i}\left[\ee^{-1/\tau}-1\right]_{\tau=t}=
\frac{1}{2\pi i}\left[\ee^{-1/\tau}\right]_{\tau=t}.
\end{align*}
In this case, the inverse Radon transform exists and is unique, 
see~\cite[Theorem 2.4]{KT95}. It can be explicitly calculated and seen
to fulfill the original differential equation. The second hyperfunction
solution can be recovered by the following trick: Changing the sign of the
$G_-$-part of the canonical defining function for 
$-2\pi i\delta(t)=[G_+(\tau_+),G_-(\tau_-)]$, where 
$G_+(\tau)=G_-(\tau)=1/\tau$, essentially gives the Cauchy 
principal value\ie the finite part distribution associated with $1/t$ and 
likewise for the higher derivatives. We therefore consider the new ansatz
\[
f_2(\omega,t)=\sum_{n=0}^\infty h_n \FP\frac{1}{t^n}
\]
(where $\FP$ denotes Hadamard's finite part, 
and with the convention $\FP 1 = 1$). Inserting
\[
t^2\frac{\dd}{\dd t}\FP\frac{1}{t^{n+1}}=-(n+1)\FP\frac{1}{t^n}
\]
yields
\[
\left(t^2\frac{\dd}{\dd t}-1\right)\sum_{n=0}^\infty h_n\FP\frac{1}{t^{n+1}}=
-h_0-\sum_{n=1}^\infty \left(h_n(n+1)+h_{n-1}\right)\FP\frac{1}{t^n}=0.
\]
With the initial condition $h_0=1$ we get
$ h_n=(-1)^n/(n+1)!$, $n=0,1,\ldots$.
On subtracting the constant hyperfunction
$1$ to compensate the term $h_0$ we get the second hyperfunction solution
by a converging series of defining functions
\[
f_2(\omega,t)=1+\sum_{n=1}^\infty \frac{(-1)^n}{n!}\FP\frac{1}{t^n}=
\frac{1}{2}\left(\ee^{-1/(t+i0)}+
\ee^{-1/(t-i0)}\right).
\]
These two hyperfunction solutions span, together with the classical
solution $f_3(t)=\ee^{-1/t}$, $t>0$, continued by $0$ to $(-\infty,0]$, 
the solution space of the original equation. 
\end{nxmpl}
\appendix
\section{Proof of the Assertions of Section~\ref{sect:ito}}
\label{app:proofs}
\subsection{Proof of Proposition~\ref{prop:DFS}}
We refer to~\cite{KOM73B} and~\cite{b:FW68} for the essentials of locally 
convex spaces defined by limits of inductive or projective sequences. 
For the corresponding results in the context of ultradistributions 
see~\cite{PET78,PET79}.
To show that \HOIDX{\PMI}{}{K} is DFS, respectively that
\HOIDX{\PMI}{}{W} is FS for every compactum $K\subset\Qn$, respectively for every open set
$W\subset\Qn$, it is enough to show the following lemma.
\begin{lemm}
The natural inclusion mapping $\rho\colon \HOBO[m]{L}\hookrightarrow\HOBO[n]{K}$ is a compact mapping whenever $m<n$ and $K\Subset L\subset\Qn$ are compact.
\end{lemm}
\begin{proof}
Let $\{f_q\}\subset\HOBO[m]{L}$ be a bounded sequence. Then, we have the estimate
\[
\sup_{z\in U} \ABS{f_q(z)}\leq M_S\NO{f_q}{n}{K}\leq N_S\NO{f_q}{m}{L},
\]
which holds, with constants $M_S$, $N_S>0$, for every compact set
$S\subset K_\Cn^\circ$, every open neighborhood $U$ of $S$ in
$K_\Cn$, and uniformly in $q$. Under these conditions,
\cite[Corollary 2.2.5]{b:HOE66} implies that there is a
subsequence $f_{q_k}$ converging to a limit
$f\in\smash{\HOIDX{}{}{K_\Cn^\circ}}$ uniformly on compact
subsets of $\smash{K_\Cn^\circ}$. The second inequality above
shows that $f_{q_k}$ converges to $f$ in \HOBO[n]{K}. This shows
the assertion.
\end{proof}
To show nuclearity of the four types of spaces in question, we
proceed as in the proof of~\cite[Proposition 2.12]{NAG81A}
or~\cite[Proposition 2.1.3]{ITO88}, by first showing nuclearity of
yet another space. Let $U\subset\Qn$ be open and define for every
integer $s$ the Fr\'{e}chet space
\begin{gather*}
  \HOIDX{\underline{s}}{}{U}\DEF
  \Bigl\{f\in\HOIDX{}{}{U_{\Cn}} \Bigm|
  \NODX{f}{s,p}{K}<\infty,\ K\subset U\text{ compact},\
  p=2,3,\ldots\Bigr\},\\
\intertext{with the following system of seminorms for $p\geq2$:}
  \NODX{f}{s,p}{K}\DEF
  \int_{K_{\Cn}} \ABS{f(z)}M_p(z)\dd\lambda(z), \quad\text{where}\quad
M_p(z)= (1+\ABS{\RE z})^{-s-\ABS{s}/p}.
\end{gather*}
Here, $\lambda$ denotes the Lebesgue measure on \Cn.
\begin{lemm}
The space \HOIDX{\underline{s}}{}{U} is nuclear.
\end{lemm}
\begin{proof}
Repeat the proof of~\cite[Proposition 2.11]{NAG81A}, or
\cite[Lemma 2.1.4]{ITO88} with minor modifications.
\end{proof}
Now, choose a fundamental system $\{U_m\}_{m\in\NN}$ of open neighborhoods of  $K$ such
that $U_{m+1}\Subset U_m$ and further an exhausting sequence $\{\smash[t]{\ICLOS{V_m}}\}_{m\in\NN}$
for $W$ such that $V_m\Subset V_{m+1}$ are open sets in $W$.
\begin{lemm}
\label{lemm:auxrep}
There are linear, topological isomorphisms
\[
    \HOPI{K}\cong\varinjlim\HOIDX{\underline{m}}{}{U_m},\quad\text{and}\quad
    \HOMI{W}\cong\varprojlim\HOIDX{\underline{-m}}{}{V_m}.
\]
\end{lemm}
\begin{proof}
We have a continuous inclusion
$\HOBO[m]{\smash[t]{\ICLOS{U_m}}}\hookrightarrow
\HOIDX{\underline{m+\AIX{n}}}{}{U_{m+\AIX{n}}}$. On the other
hand, an application of Cauchy's integral formula as in the proof
of~\cite[Lemma 2.1.4]{ITO88} yields an embedding
$\HOIDX{\underline{m}}{}{U_m}\hookrightarrow\HOBO[m+1]{\smash[t]{\ICLOS{U_{m+1}}}}$.
This suffices to show equivalence of the inductive limits in the
tempered case. The asymptotic case follows by similar
considerations.
\end{proof}
The permanence properties of nuclearity, see~\cite[Proposition
50.1]{b:TRE67}, then imply all assertions of
Proposition~\ref{prop:DFS}.
\subsection{Proof of Proposition~\ref{prop:tensor}}
Since the proof is mainly an application of general facts about
locally convex spaces and their tensor products, we refer
the reader to the proofs of~\cite[Proposition
2.1.7--2.1.10]{ITO88}, which in turn follows~\cite[Proposition
3.6]{NAG81A}, and merely note the only modification that has to be
inserted: There is an auxiliary space $\mathcal{E}_\ast$ of
$C^\infty$-functions which has to be replaced by
\begin{align*}
\mathcal{E}_{-\infty}(U) &\DEF \Bigl\{f\in C^\infty(U_\Cn)\Bigm|
  \forall\gamma>0\colon\sup_{z\in
  K_\Cn}\ABS{\smash[t]{f^{(\alpha)}(z)(1+\ABS{z})^\gamma}}<\infty \Bigr\},\\
\mathcal{E}_\infty(U) &\DEF \Bigl\{f\in C^\infty(U_\Cn)\Bigm|
  \exists\gamma\in\RR\colon\sup_{z\in
  K_\Cn}\ABS{\smash[t]{f^{(\alpha)}(z)(1+\ABS{z})^\gamma}}<\infty
\Bigr\},
\end{align*}
for the asymptotic and tempered case of the assertions
respectively, and any open subset $U$ of \Qn. The conditions in
the definitions above are meant to hold for every compact set
$K\subset U$ and every real partial derivative $f^{(\alpha)}$,
$\alpha\in\NN_0^{2\AIX{n}}$. The proof can then be carried out as
in the reference indicated above, using further the
results~\cite[Theorem 39.2 and Proposition 36.1]{b:TRE67},
\cite[Corollary 1 of Lemma A]{ITO79}, and~\cite[Theorem 5 of
§4]{BM73}.
\subsection{Proof of Theorem~\ref{theo:H1vanish}}
This proof follows~\cite[Theorem 2.1.14]{ITO88}. It heavily
depends on a notion originally designed for slowly increasing\ie
infra-ex­po­nen­ti­al functions in the context of Fourier
hyperfunctions, see~\cite{KAW70}, but which is also applicable in
our case due to its mainly geometric nature.
\begin{defi}[{\cite[Definition 5.1]{NAG81A}}]
An open subset $V$ of \Qn is called an \textbf{\HOIE{}-pseu­do­con­vex} set if
\begin{enumerate}
\item $\sup\bigl\{|\IM z''|, |\IM z'''|-|\RE z'''|\bigm| z\in V_\Cn\bigr\}<\infty$.
\item There exists a smooth, plurisubharmonic function $\varphi(z)$ on $V_\Cn$ which is bounded on $L_\Cn$ for every compactum $L\subset V$, and such that $\smash[t]{\ICLOS{V_c}}$ is compact in $V$, where $V_c\DEF\bigl\{z\in v_\Cn\bigm|\varphi(z)<c\bigr\}$.
\end{enumerate}
\end{defi}
Condition i) only states that a pseudoconvex domain is of finite
width around \Dn in the second part of the coordinates
$\IX{n}=(n_1,n_2,n_3)$, and of finite inclination above \Dn in
the third part. The second condition is a direct generalization
of the notion of pseudoconvexity in the complex domain,
see~\cite[Theorem 2.6.7]{b:HOE66}. The following is an analogue of
Grauert's Theorem:
\begin{theo}[{\cite[Theorem 5.3]{NAG81A}, \cite[Theorem 2.1.13]{ITO88}}]
For every open\hfill\linebreak $S\subset\Dn$ and an open neighborhood $U\subset\Qn$ of $S$, there exists an \HOIE{}-pseudoconvex set $V\subset U$ such that $S=V\cap\Dn$.
\end{theo}
Thus the compact set $K$ in the statement of Theorem~\ref{theo:H1vanish} has a fundamental system of neighborhoods in \Qn consisting of \HOIE{}-pseudoconvex sets. Therefore, we have only to prove that $\RELCO{1}{U}{\Qn}{\HOPMI{}}=0$ for every \HOIE{}-pseudoconvex set $U\subset\Qn$, such that $\ABS{\IM{z'}}<\EPS$, $\ABS{\IM{z''}}<\EPS$, and $\ABS{\IM{z'''}}<\EPS(1+\ABS{\RE{z'''}})$ on $U_\Cn$ for some \EPS sufficiently small for our later purposes. Since $U$ is paracompact, the relative cohomology groups coincide with the \v{C}ech cohomology groups and it remains only to show $\RELCO{1}{}{\{U_j\}_{j\in\NN}}{\HOPMI{}}=0$ for any locally finite covering $\{U_j\}_{j\in\NN}$ of $U$ such that $V_j=U_{j\Cn}$ is \HOIE{}-pseudoconvex.

We prove the assertion for \SHOPI\ only, since a very similar argument applies in the case of asymptotic functions. Let $C^s(Z^{\text{loc}}(\{V_j\}))$ be the space of cochains $c=\{c_J\mid J=(j_0,\ldots,j_s)\in\NN^{s+1}\}$ such that
\begin{enumerate}
\item $\overline{\partial}c_J=0$ on $V_J=V_{j_0}\cap\ldots\cap V_{j_s}$,
\item $\smash[b]{\sum_{J\in M}\int_{V_J}}|c_J|^2\dd\lambda<\infty$, for any finite
subset $M$ of $\NN^{s+1}$.
\end{enumerate}
If $d=\{d_{ij}\}$ represents a cocycle in $\RELCO{1}{}{\{U_j\}_{j\in\NN}}{\HOPI{}}$ then there is a $\gamma>0$ such that $d_{ij}\cdot j_{-\gamma}|_{V_{ij}}$ is in $C^1(Z^{\text{loc}}(\{V_j\}))$, where the multiplier $j_\gamma$ is defined for $\gamma\in\RR$ by
$
j_\gamma\DEF\bigl( 1+{z''}^2+{z'''}^2  \bigr)^{\gamma/2}.
$
Here, we assume that the \EPS above is small enough such that $j_{-\gamma}$ and $j_{-\gamma-2}$ (to be used below) are holomorphic on $U_\Cn$\ie we implicitly consider the simultaneous inductive limit over neighborhoods $U$ of $K$ and increasing growth order of forms in \RELCO{1}{}{U}{\HOPI{}}.
Define $c=\{c_{ij}\}\in C^1(Z^{\text{loc}}(\{V_j\}))$ such that $\delta c=0$, by $c_{ij}\DEF d_{ij}\cdot j_{-\gamma-2}|_{V_{ij}}$, where $\delta$ is the coboundary operator. Let $\{\chi_j\}$ be a smooth partition of unity subordinate to $\{V_j\}$, and set $b_j\DEF\sum_i\chi_i c_{ij}$. Then $\delta c=0$ implies $\delta b=c$ and consequently $\delta\dbar b=\dbar c=0 $ by condition i).
Since $\sum\chi_i=1$ and $\chi_i\geq 0$ we have
\[
\int_{V_j}|b_j|^2\dd\lambda\leq
\sum_i\int_{V_i}\chi_i|c_{ij}|^2\dd\lambda<\infty
\]
by the triangle inequality and condition ii). Since $U$ is \HOIE{}-pseudoconvex, there exists a smooth plurisubharmonic function $\psi$ on $V\DEF U_\Cn$ that satisfies the two conditions (1) $\sum_j|\dbar\chi_j(z)|\leq\exp(\psi(z))$ and (2) $\sup_{K_\Cn}\psi\leq C_K$ for every relatively compact subset $K$ of $U$. Thus, (1) and condition ii) above, together with the definition of $c$ and $b$ imply
\[
\sum_{j\in M} \int_{V_j} |\dbar b_j|^2(1+|z|^2)^2\exp(-\psi(z))\dd\lambda<\infty,
\]
for all finite sets $M\subset\NN$. Since $\dbar b$ is closed, it
defines a global section $f$ on $U_\Cn$. This and the last
estimate allow us to make use of~\cite[Theorem 4.4.2]{b:HOE66},
and conclude that there exists a smooth section $u$ on $U_\Cn$
such that $\dbar u=f$, and with $ \int_{K_\Cn}
|u|^2\dd\lambda<\infty $ for all relatively compact subsets $K$
of $U$ (this theorem increases the growth order by two, but
remember the estimate on the integral over $|\dbar b_j|$ above.)
Set $c_j'\DEF b_j-u|_{V_j}$. Then $\dbar c_j'=0$, $\delta
c'=\delta b=c$, and $c'$ is an element of
$C^1(Z^{\text{loc}}(\{V_j\}))$. Finally, with $d_j'\DEF c_j'\cdot
j_{\gamma+2}|_{V_j}$ we find that the collection $d'\DEF\{d_j'\}$
is a subset of \HOPI{} such that $\delta d'=d$, which means that
$d=0$ in $\RELCO{1}{}{\{U_j\}_{j\in\NN}}{\HOPI{}}$. Since $d$ was
a generic element, we find
$\RELCO{1}{}{\{U_j\}_{j\in\NN}}{\HOPI{}}=0$.
\subsection{Proof of Theorem~\ref{theo:runge}}
The proof follows the original one in the case of Fourier
hyperfunctions~\cite[Theorem 2.2.1]{KAW70}, see also~\cite[proof
of Theorem 3.1]{NAG81A} and~\cite[Section 2.2]{ITO88}.

Again, we consider only the case \SPPI{}. Since \Dn is
$\sigma$-compact, the assertion of the theorem factorizes in the
first and the last two variables and it is sufficient to prove:
\begin{theo}
\label{theo:densereduct}
If $L=\{x'\in\RR^{n_1}\mid |x'|\leq a\}×\DD^{n_2}×\DD^{n_3}$ contains
the compact set $K\subset\Dn$ for some $a>0$ then \SPPI{L} is dense in \SPPI{K}.
\end{theo}
\newcommand{\OL}[2]{\HOIDX{#1}{2,\text{loc}}{#2}}
\newcommand{\LL}[2]{\SetConstructor{L}{#1}{2,\text{loc}}{#2}}
\newcommand{\LC}[2]{\SetConstructor{L}{#1}{2,\text{c}}{#2}}
For $W\subset\Qn$ open and $\eta>0$, we define \OL{\eta}{W} to be
the space of all holomorphic functions on $W_\Cn$ such that for
all compacta $K\subset W$ holds: $
\int_{K_\Cn}|f|^2(1+|z''|+|z'''|)^{-\eta}\dd\lambda<\infty $.
This space is a FS-space (consider any exhaustive sequence of
compacta for $W$), and we have
\begin{lemm}
If $\{W_j\}$ is a fundamental system of neighborhoods for a compact subset
$K$ in \Qn, then there is a linear topological isomorphism
\[\HOPI{K}\cong\varinjlim \OL{j}{W_j}.\]
\end{lemm}
To see this, note that for sufficiently large $j$, there are numbers $k$, $l$ such that
there exist continuous inclusions
$\HOBO[j]{W_j}\hookrightarrow\OL{j+k}{W_{j+k}}$ and
$\OL{j+k}{W_{j+k}}\hookrightarrow\HOBO[j+k+l]{W_{j+k+l}}$.

Now, let $\{W_j\}$ and $\{V_j\}$ be fundamental systems of
neighborhoods for the compact sets $L\subset K\subset\Dn$ of
Theorem~\ref{theo:densereduct} respectively, such that
$V_j\subset W_j$ for all $j$. Then
$\HOPI{L}\cong\varinjlim\OL{j}{W_j}$ and
$\HOPI{K}\cong\varinjlim\OL{j}{V_j}$. We are done if we are able
to show that \OL{j}{W_j} is dense in \OL{j}{V_j} for sufficiently
large $j$, cf.~\cite[Lemma 2.2.7]{ITO88}. For this, it trivially
suffices to show that \OL{l}{W_j} is dense in \OL{j}{V_j} for
$l\leq j$ sufficiently large. Setting for brevity $W=W_j$,
$V=V_j$, the Hahn--Banach Theorem tells us that we only need to
prove the following: If $\mu\in\OL{j}{V}'$ and $\SPROD{\mu}{v}=0$
for all $v\in\OL{l}{W}$ then $\mu=0$.

The solution of our problem requires finding solutions of the
dual Cauchy-Riemann differential equation with growth conditions.
These solutions traditionally live in $L^2$--spaces, and so we
have to make another definition: \LL{\eta}{W} is the space of
locally square integrable functions on $W_\Cn$ that satisfy the
same integrability condition as the functions in \OL{\eta}{W}.
\LL{\eta}{W} is a FS*--space and \OL{\eta}{W} is a closed
subspace of it. The dual space of \LL{\eta}{W} is the space
\LC{-\eta}{W} of functions $f\in\LL{}{W_\Cn}$ with compact
support in $W$, and such that
$\int_{W_\Cn}|f|^2(1+|z''|+|z'''|)^{\eta}\dd\lambda<\infty$.

Now, since \OL{j}{V} is a closed subspace of \LL{j}{V}, there exists a
re­pre­sen­ta­tive $u\in\LC{-j}{V}$ of $\mu\in\OL{j}{V}'$\ie
\[
\SPROD{\mu}{v}=\int_{V_\Cn} v\overline{u}\dd\lambda,\quad
\text{for } v\in\OL{j}{V}.
\]
Thus we only have to prove that $u$ is orthogonal to \OL{j}{V}.

Set $T=\ICLOS{\supp(u)}$. By~\cite[Lemma 2.2.6]{ITO88}, there
exists a neighborhood $U$ of $T$ which is relatively compact in
$V$, and a strictly plurisubharmonic smooth function $\theta$ on
$W_\Cn$ such that i) $\theta(z)<0$ on $T_\Cn$, ii) $\theta(z)>0$
on a neighborhood $N$ of $(\partial U)_\Cn$, and iii)
$\sup_{L_\Cn}\theta(z)<\infty$ for any relatively compact set
$L\subset W$. Then, \cite[Theorem 2.3.2]{HOE65} ensures that there
exists a form $f\in L_2^{(0,1)}(W_\Cn;(2-j)\log(1+\ABS{z}))$ such
that $u=\vartheta f$, where $\vartheta$ is the adjoint of the
$\partial$-operator. Furthermore, $\supp(f)\subset\{z\in
W_\Cn\mid \theta(z)\leq 0\}$. Choose a smooth function $\chi$ on
$W_\Cn$ such that $0\leq\chi(z)\leq1$, $\chi(z)=1$ on $T_\Cn$,
$\chi(z)=0$ on $\smash[t]{\overline{(U\cup N)_\Cn}}$,
$\supp(\overline{\partial}\chi)\subset N$, and
$\sup\ABS{\smash{\overline{\partial}\chi}}<\infty$. Then, for
every $w\in\OL{j-2}{V}$ holds
\begin{align*}
\SPROD{\mu}{w}&=
\int_{V_\Cn} w \overline{u}\dd\lambda =
\int_{W_\Cn} (\chi w) \overline{u}\dd\lambda\\
& = \int_{W_\Cn} (\chi w) \overline{\vartheta f}\dd\lambda =
 \int_{W_\Cn} \overline{\partial} (\chi w) \overline{f}\dd\lambda=0,
\end{align*}
by the constraints on $f$ and $\chi$. A trivial analogue of
of~\cite[Lemma 2.2.5]{ITO88} implies that \OL{j-2}{V} is dense in
\OL{j}{V} if we choose the shape of $V$ as in~\cite[p.\
234]{ITO88}. This completes the proof of
Theorem~\ref{theo:densereduct}.
\newcommand{\noopsort}[1]{} \newcommand{\singleletter}[1]{#1}
\providecommand{\bysame}{\leavevmode\hbox to3em{\hrulefill}\thinspace}
\providecommand{\MR}{\relax\ifhmode\unskip\space\fi MR }
\providecommand{\MRhref}[2]{%
  \href{http://www.ams.org/mathscinet-getitem?mr=#1}{#2}
}

\begin{thebibliography}{10}

\bibitem{BM73}
Klaus-Dieter Bierstedt and Reinhold Meise, \emph{Lokalkonvexe {U}nterr{ä}ume
  in topologischen {V}ektorr{ä}umen und das {$\epsilon$}-{P}rodukt},
  Manuscripta Math. \textbf{8} (1973), 143--172.

\bibitem{BS74}
E.~Br{ü}ning and P.~Stichel, \emph{On the {E}quivalence of {S}caling,
  {L}ight-cone {S}ingularities, and {A}symptotic {B}ehavior of the
  {J}ost-{L}ehmann {S}pectral {F}unction}, Comm.\ Math.\
  Phys.\ \textbf{36} (1974), 137--156.

\bibitem{BS75}
\bysame, \emph{Asymptotics and {L}ight-cone {S}ingularities in {Q}uantum
  {F}ield {T}heory}, Proceedings of the {I}nternational {S}ymposium on
  {M}athematical {P}roblems in {T}heoretical {P}hysics (Kyoto, Japan, January
  19--23, 1975) (H.~Araki, ed.), Springer-{V}erlag, 1975, pp.~72--84.

\bibitem{BN89}
Erwin Br{ü}ning and Shigeaki Nagamachi, \emph{Hyperfunction {Q}uantum {F}ield
  {T}heory: {B}asic {S}tructural {R}esults}, J.~Math.\ Phys.\
   \textbf{30} (1989), no.~10, 2340--2359.

\bibitem{BN98}
\bysame, \emph{Closure of {F}ield {O}perators, {A}symptotic {A}belianess, and
  {V}acuum {S}tructure in {H}yperfunction {Q}uantum {F}ield {T}heory}, J.~Math.\ 
  Phys.\ \textbf{39} (1998), 5098--5111.

\bibitem{CKL95}
Soon-Yeong Chung, Dohan Kim, and Eun~Gu Lee, \emph{Schwartz {K}ernel {T}heorem
  for the {F}ourier {H}yperfunctions}, Tsukuba J.\ Math.\ \textbf{19} (1995),
  no.~2, 377--385.

\bibitem{CIO90}
Ioana Cioranescu, \emph{Moment sequences for ultradistributions}, Math.~Z.\
\textbf{204} (1990), 391--400.

\bibitem{CT74}
Florin Constantinescu and Willi Thaleimer, \emph{Euclidean {G}reen´s
  {F}unctions for {J}affe {F}ields}, Comm.\ Math.\ Phys.\
  \textbf{38} (1974), no.~2, 299--316.

\bibitem{CT79}
\bysame, \emph{Ultradistributions and {Q}uantum {F}ields: {F}ourier-{L}aplace
  {T}ransforms and {B}oundary {V}alues of {A}nalytic {F}unctions}, Rep.\
  Math.\ Phys.\ \textbf{16} (1979), no.~2, 167--180.

\bibitem{DE94}
A.~L. Dur{á}n and R.~Estrada, \emph{Strong {M}oment {P}roblems for {R}apidly
  {D}ecreasing {S}mooth {F}unctions}, Proc.\ Amer.\ Math.\ Soc.\
  Society \textbf{120} (1994), no.~2, 529--534.

\bibitem{EGV98}
R.~Estrada, J.~M. Gracia-Bond{í}a, and J.~C. V{á}rilly, \emph{On
  {S}ummability of {D}istributions and {S}pectral {G}eometry}, Comm.\
  Math.\ Phys.\ \textbf{191} (1998), 219--248.

\bibitem{EK90}
R.~Estrada and R.~P. Kanwal, \emph{A {D}istributional {T}heory for {A}symptotic
  {E}xpansions}, Proc.\ Roy.\ Soc.\ Lond.\ Ser.~A \textbf{428} (1990), 399--430.

\bibitem{b:EK94}
R.~Estrada and R.~P. Kanwal, \emph{Asymptotic {A}nalysis: {A} {D}istributional
  {A}pproach}, Birkh{ä}user--{V}erlag, Boston, Basel, Stuttgart, 1994.

\bibitem{b:FW68}
Klaus Floret and Joseph Wloka, \emph{Einf{ü}hrung in die {T}heorie der
  lokalkonvexen {R}{ä}ume}, Springer-{V}erlag, Berlin, {H}eidelberg, {N}ew
  {Y}ork, 1968.

\bibitem{b:GS64}
I.~M. Gelfand and G.~E. Shilov, \emph{Generalized {F}unctions}, vol.~2,
  Academic {P}ress, London, {S}an {D}iego, 1964.

\bibitem{GLS68}
A.~Grossman, G.~Loupias, and E.~M. Stein, \emph{An {A}lgebra of
  {P}seudodifferential {O}perators and {Q}uantum {M}echanics in {P}hase
  {S}pace}, Ann. Inst. Fourier \textbf{18} (1968), 343--368.

\bibitem{HAS61}
Morisuke Hasumi, \emph{Note on the {$n$}-{D}imensional {T}empered
  {U}ltra-{D}is­tri­bu­tions}, T\^{o}hoku Math.~J.\ \textbf{13} (1961),
  94--104.

\bibitem{b:HEL84}
Sigurdur Helgason, \emph{Groups and Geometric Analysis. Integral Geometry, 
Invariant Differential Operators, and Spherical Functions.}, 
Pure and Applied Mathematics, vol. 113, Academic {P}ress,
Orlando, Florida, 1984.

\bibitem{HOE65}
Lars Hörmander, \emph{{$L^2$}-estimates and {E}xistence {T}heorems for the
  $\overline{\partial}$--{O}perator}, Acta Math.\ \textbf{113} (1965),
  89--152.

\bibitem{b:HOE66}
Lars H{ö}rmander, \emph{An {I}ntroduction to {C}omplex {A}nalysis in
  {S}everal {V}ariables}, Van {N}ostrand {R}einhold {C}ompany, {N}ew {Y}ork,
  {T}oronto, {L}ondon, {M}elbourne, 1966.

\bibitem{ITO79}
Yoshifumi Ito, \emph{On the {T}heory of {V}ector {V}alued
  {H}yperfunctions}, J.~Math.\ Tokushima Univ.\ \textbf{13} (1979), 29--51.

\bibitem{ITO88}
\bysame, \emph{Fourier {H}yperfunctions of {G}eneral {T}ype}, J.~Math.\ Kyoto
  Univ.\ \textbf{28} (1988), 213--265.

\bibitem{ITO92}
\bysame, \emph{Vector {V}alued {F}ourier {H}yperfunctions}, J.~Math.\ Kyoto
  Univ.\ \textbf{32} (1992), 259--285.

\bibitem{JS88}
L.~Jakobczyk and F.~Strocchi, \emph{Euclidean {F}ormulation of {Q}uantum
  {F}ield {T}heory {W}ithout {P}ositivity}, Comm.\ Math.\
  Phys.\ \textbf{119} (1988), 529--541.

\bibitem{JUN78}
K.~Junker, \emph{Vektorwertige {F}ourierhyperfunktionen}, Diplomarbeit,
  Universit{ä}t D{ü}sseldorf, 1978.

\bibitem{b:KAN88}
Akira Kaneko, \emph{Introduction to {H}yperfunctions}, Kluwer {A}cademic
  {P}ublishers, {D}ordrecht, 1988.

\bibitem{KT95}
Akira Kaneko and Takashi Takiguchi, \emph{Radon {T}ransform of {H}yperfunctions
  and {S}upport {T}heorem}, Hokkaido Math.~J.\ \textbf{24} (1995),
  63--103.

\bibitem{KK80}
Masaki Kashiwara and Takahiro Kawai, \emph{Second microclocalization and
  asymptotic expansions}, Complex Analysis, Microlocal Calculus and
  Relativistic Quantum Theory  (D.~Iagolnitzer, ed.), Lecture
  {N}otes in {P}hysics, vol. 126, Proceedings of the colloquium at the Centre
  de Physique, Les Houches, September 1979, Springer-{V}erlag, 1980, pp.~21--76.

\bibitem{b:KKK86}
Masaki Kashiwara, Takahiro Kawai, and Tatsuo Kimura, \emph{Foundations of
  {A}lgebraic {A}nalysis}, Princeton {U}niversity {P}ress, Princeton, {N}ew
  {J}ersey, 1986.

\bibitem{KAW70}
Takahiro Kawai, \emph{On the {T}heory of {F}ourier {H}yperfunctions and its
  {A}pplications to {P}artial {D}ifferential {E}quations with
  {C}onstant {C}oefficients}, J.~Fac.\ Sci.\ Univ.\ Tokyo Sect.~IA Math.\ 
  \textbf{17} (1970), 467--515.

\bibitem{KOM73B}
Hikosaburo Komatsu, \emph{Ultradistributions, {I} {S}tructure {T}heorems and a
  {C}haracterization}, J.~Fac.\ Sci.\ Univ.\ Tokyo Sect.~IA Math.\ \textbf{20}
  (1973), 25--105.

\bibitem{KOM89}
\bysame, \emph{Microlocal {A}nalysis in {G}evrey {C}lasses and in {C}omplex
  {D}omains}, Microlocal {A}nalysis and {A}pplications  
  (Lamberto Cattabriga and Luigi Rodino, eds.), Lecture {N}otes in
  {M}athematics, vol. 1495, Springer-{V}erlag, Berlin,
  pp.~161--236, 2nd session of the {C}entro {I}nternazionale
  {M}atematico {E}stivo ({C.I.M.E.}), Montecattini Terme,
  Italy.

\bibitem{b:MOR93}
Mitsuo Morimoto, \emph{An {I}ntroduction to {S}ato's {H}yperfunctions},
  American {M}athematical {S}ociety, Providence, {R}hode {I}sland, 1993.

\bibitem{MS92}
U.~Moschella and F.~Strocchi, \emph{The {C}hoice of {T}est {F}unctions in
  {G}auge {Q}uantum {F}ield {T}heory}, Lett.\ Math.\ Phys.\
  \textbf{24} (1992), 103--113.

\bibitem{NM76A}
S.~Nagamachi and N.~Mugibayashi, \emph{Hyperfunction {Q}uantum {F}ield
  {T}heory}, Comm.\ Math.\ Phys.\ \textbf{46} (1976),
  119--134.

\bibitem{NM76B}
\bysame, \emph{Hyperfunction {Q}uantum {F}ield {T}heory {II}. Euclidean Green's Functions}, 
Comm.\ Math.\ Phys.\ 
 \textbf{49} (1976), 257--275.

\bibitem{NM77}
\bysame, \emph{Quantum {F}ield {T}heory in {T}erms of {F}ourier
  {H}yperfunctions}, Publ.\ Res.\ Inst.\ Math.\ Sci.\ 
  Kyoto Univ.\ \textbf{12} (1977), 309--341,
  Supplement.

\bibitem{NM79}
\bysame, \emph{The {H}aag--{R}uelle {F}ormulation of {S}cattering in
  {H}yperfunction {Q}uantum {F}ield {T}heory}, Rep.\ Math.\
  Phys.\ \textbf{16} (1979), 181--201.

\bibitem{NM86}
\bysame, \emph{Hyperfunctions and {R}enormalization}, J.~Math.\
  Phys.\ \textbf{27} (1986), no.~3, 832--838.

\bibitem{NAG81A}
Shigeaki Nagamachi, \emph{The {T}heory of {V}ector {V}alued {F}ourier
  {H}yperfunctions of {M}ixed {T}ype~{I}}, Publ.\ Res.\ Inst.\ Math.\ Sci.\ 
  Kyoto Univ.\ \textbf{17}
  (1981), 25--63.

\bibitem{NAG81B}
\bysame, \emph{The {T}heory of {V}ector {V}alued {F}ourier {H}yperfunctions of
  {M}ixed {T}ype~{II}}, Publ.\ Res.\ Inst.\ Math.\ Sci.\ 
  Kyoto Univ.\ \textbf{17} (1981), 65--93.

\bibitem{b:OLV97}
Frank~W.~J. Olver, \emph{Asymptotics and {S}pecial {F}unctions}, Academic
  {P}ress, London, {S}an {D}iego, 1997.

\bibitem{PM73}
Young~Sik Park and Mitsuo Morimoto, \emph{Fourier {U}ltra-{H}yperfunctions in
  the {E}uclidean {$n$}-{S}pace}, J.~Fac.\ Sci.\ Univ.\ Tokyo Sect.~IA Math.\
  \textbf{20} (1973), 121--127.

\bibitem{PET78}
Hans--Joachim Petzsche, \emph{Die {N}uklearität der {U}ltradistributionsräume 
und der {S}atz vom {K}ern. {I}}, Manuscripta Math.\ 
\textbf{24} (1978), no.~2, 133--171.

\bibitem{PET79}
\bysame, \emph{Die {N}uklearität der {U}ltradistributionsräume 
und der {S}atz vom {K}ern. {II}}, Manuscripta Math.\ 
\textbf{27} (1979), no.~3, 221--251.

\bibitem{PET84}
\bysame, \emph{Generalized functions and the boundary values 
of holomorphic functions}, J.~Fac.\ Sci.\ Univ.\ Tokyo Sect.~IA Math.\ 
\textbf{31} (1984), no.~2, 391--431.

\bibitem{PIE88}
D.~Pierotti, \emph{The {E}xponential of the {T}wo-{D}imensional {M}assless
  {S}calar {F}ield as an {I}nfrared {J}affe {F}ield}, Lett.\ Math.\
  Phys.\ \textbf{15} (1988), 219--230.

\bibitem{SAB85}
Yutaka Saburi, \emph{Fundamental {P}roperties of {M}odified {F}ourier
  {H}yperfunctions}, Tokyo J.\ Math.\ 
  \textbf{8} (1985), no.~1, 231--273.

\bibitem{SAT59}
Mikio Sato, \emph{Theory of {H}yperfunctions, {I}}, 
J.~Fac.\ Sci.\ Univ.\ Tokyo Sect.~IA Math.\ 
\textbf{8} (1959), no.~1, 139--193.

\bibitem{SAT60}
\bysame, \emph{Theory of {H}yperfunctions, {II}}, 
J.~Fac.\ Sci.\ Univ.\ Tokyo Sect.~IA Math.\ 
\textbf{8} (1960), no.~2, 387--437.

\bibitem{SCH97B}
Andreas~U. Schmidt, \emph{Euclidean {R}econstruction in {Q}uantum {F}ield
  {T}heory: {B}etween {T}empered {D}istributions and {F}ourier
  {H}yperfunctions}, Preprint, Johann Wolfgang Goethe-Universit{ä}t,
  Frankfurt am Main, Germany, 1997, 
  \href{http://xxx.lanl.gov/abs/math-ph/9811002}{\texttt{math-ph/9811002}}.

\bibitem{SCH97A}
\bysame, \emph{Mathematical {P}roblems of {G}auge {Q}uantum {F}ield {T}heory:
  {A} {S}urvey of the {S}chwinger {M}odel}, Universitatis Iagellonicae Acta
  Mathematica \textbf{Fasciculus XXXIV} (1997), 113--134,
  \href{http://xxx.lanl.gov/abs/hep-th/9707166}{\texttt{hep-th/9707166}}.

\bibitem{b:SCH99}
\bysame, \emph{Asymptotische {H}yperfunktionen, temperierte {H}yperfunktionen
  und asymptotische {E}ntwicklungen}, Logos-Verlag, Berlin, 1999, Dissertation,
  Frankfurt am Main, 1999.

\bibitem{SCH02}
\bysame, \emph{{I}nfinite infrared regularization and a state space for the 
{H}eisenberg algebra}, J.~Math.\ Phys.\ \textbf{43} (2002), no.~1, 243--259,
\href{http://xxx.lanl.gov/abs/math-ph/0105042}{\texttt{math-ph/0105042}}.

\bibitem{SCH03}
\bysame, \emph{Phragmén--Lindelöf Principles of type $L^p$},
Appl.\ Math.\ E-Notes \textbf{3} (2003), 178--182,
\href{http://www.math.nthu.edu.tw/~amen/2003/020923-4.pdf}{http://www.math.nthu.edu.tw/~amen/2003/020923-4.pdf}

\bibitem{SEB58}
J.~{Sebasti{ã}o}~e Silva, \emph{Les fonctions analytiques commes
  ultradistributiones dans le calcul operationn\'{e}l}, Math.\ Ann.\
  \textbf{136} (1958), 58--96.

\bibitem{SOL97A}
M.~A. Soloviev, \emph{An {E}xtension of {D}istribution {T}heory and of the
  {P}ayley-{W}iener-{S}chwartz {T}heorem {R}elated to {Q}uantum {G}auge
  {T}heory}, Comm.\ Math.\ Phys.\ \textbf{184} (1997),
  579--596.

\bibitem{SOL97B}
\bysame, \emph{{W}ick-{O}rdered {E}ntire {F}unctions of the {I}ndefinite
  {M}etric {F}ree {F}ield}, Lett.\ Math.\ Phys.\ \textbf{41}
  (1997), 265--277.

\bibitem{b:STR93}
F.~Strocchi, \emph{Selected {T}opics on the {G}eneral {P}roperties of {Q}uantum
  {F}ield {T}heory}, Lecture {N}otes in {P}hysics, vol.~51, World {S}cientific Publishing,
  New Jersey, 1993.

\bibitem{b:HOO94}
G.~{'t Hooft}, \emph{Under the {S}pell of the {G}auge {P}rinciple}, World
  {S}cientific Publishing, New Jersey, 1994.

\bibitem{b:TRE67}
Fran\c{c}ois Treves, \emph{Topological {V}ector {S}paces, {D}istributions and
  {K}ernels}, Academic {P}ress, London, {S}an {D}iego, 1967.

\bibitem{b:VDZ90}
V.~S. Vladimirov, Yu.~N. Drozzinov, and B.~I. Zavialov, \emph{Tauberian
  {T}heorems for {G}eneralized {F}unctions}, Kluwer {A}cademic {P}ublishers,
  {D}ordrecht, 1988.

\bibitem{b:WAS65}
Wolfgang Wasow, \emph{Asymptotic expansions for ordinary differential
  equations}, Interscience Publishers John Wiley \& Sons, London, 1965.

\bibitem{b:WW52}
E.~T. Whittaker and G.~N. Watson, \emph{A {C}ourse of {M}odern {A}nalysis},
  Cambridge {U}niversity {P}ress, Cambridge, {UK}, 1952.

\bibitem{WIG81}
Arthur~S. Wightman, \emph{The {C}hoice of {T}estfunctions in {Q}uantum {F}ield
  {T}heory}, J.~Math.\ Anal.\ Appl.\
  \textbf{7B} (1981), 769--791, Supplement.

\bibitem{CKY98}
Yongjin Yeom, Soon-Yeong Chung, and Dohan Kim, \emph{Moment {P}roblem for
  {H}yperfunctions and {H}eat {K}ernel {M}ethod}, Preprint RIM-GARC 98--8,
  Seoul National University, Seoul, Korea, April 1998.

\end{thebibliography}
\end{document}